\documentclass[11pt]{article}

\usepackage{amstext,amssymb,amsmath,amsbsy}
\usepackage{xcolor}
\usepackage{hyperref}
\usepackage{amscd}
\usepackage{amsfonts}
\usepackage{indentfirst}
\usepackage{verbatim}
\usepackage{amsmath}
\usepackage{amsthm}
\usepackage{enumerate}
\usepackage{graphicx}
\usepackage{color, soul}
\usepackage[OT1]{fontenc}
\usepackage[latin1]{inputenc}
\usepackage[english]{babel}
\usepackage{amssymb}
\usepackage{subfig}
\usepackage{algorithm}
\usepackage{algpseudocode}

\usepackage{mathtools}

\newcommand{\R}{\mathbb{R}}
\newcommand{\N}{\mathbb{N}}
\newcommand{\ds}{\displaystyle}

\newcommand{\x}{{\bf x}}

\newcommand{\bn}{{\bf n}}

\newtheorem{remark}{Remark}[section]

\newtheorem*{Assumption*}{Assumption}

\newtheorem{problem}{Problem}[section]
\newtheorem*{problem*}{Problem}
\numberwithin{equation}{section}

\begin{document}

\title{Numerical differentiation by  the polynomial-exponential basis}

\author{
Phuong M. Nguyen\thanks{Department of Mathematics and Statistics, University of North Carolina at
Charlotte, Charlotte, NC, 28223, USA, 
\texttt{pnguye45@uncc.edu}.}
\and
Thuy T. Le\thanks{Department of Mathematics and Statistics, University of North Carolina at
Charlotte, Charlotte, NC, 28223, USA, 
\texttt{tle55@uncc.edu}.}
\and
Loc H. Nguyen\thanks{Department of Mathematics and Statistics, University of North Carolina at
Charlotte, Charlotte, NC, 28223, USA, 
\texttt{loc.nguyen@uncc.edu}, corresponding author.}
\and  Michael V. Klibanov\thanks{Department of Mathematics and Statistics, University of North Carolina at
Charlotte, Charlotte, NC, 28223, USA, 
\texttt{mklibanv@uncc.edu}.}
}


\date{}
\maketitle
\begin{abstract}
 Our objective is to calculate the derivatives of data corrupted by noise. This is a challenging task as even small amounts of noise can result in significant errors in the computation. This is mainly due to the randomness of the noise, which can result in high-frequency fluctuations. To overcome this challenge, we suggest an approach that involves approximating the data by eliminating high-frequency terms from the Fourier expansion of the given data with respect to the polynomial-exponential basis.
 This truncation method helps to regularize the issue, while the use of the polynomial-exponential basis ensures accuracy in the computation. We demonstrate the effectiveness of our approach through numerical examples in one and two dimensions.
\end{abstract}



\section{Introduction} \label{sec intr}

We revisit the problem {\it how to differentiate noisy data}. 
The problem is formulated as follows:
\begin{problem}
	Let $\Omega = (-R, R)^d$, $d \geq 1$, be an open and bounded domain of $\R^d$ where $R$ is a positive number.
	Let $f^*: \Omega \to \R$ be a function.
	Given the noisy data
	\begin{equation}
	 	f^{\delta}(\x) = f^*(\x) +  \delta\eta(\x)
		\quad \x \in \Omega,
		\label{data}
	 \end{equation} where $\eta$ is a function taking random numbers in the range $[-1, 1]$ and $\delta \in (0, 1)$ is the noise level, find an approximation of the derivatives of $f^*$. More precisely, compute
	\begin{enumerate}
	 \item the gradient vector $\nabla f^* = (\partial_{x_i} f^*)_{i = 1}^d$, 
	 \item the Hessian matrix $D^2 f^* = (\partial_{x_i x_j}^2 f^*)_{i, j = 1}^d.$
	 \end{enumerate}
	 \label{p}
\end{problem}

The task of differentiating noisy data holds great importance in various fields, including applied mathematics, statistics, data science, physics, engineering, etc. 
However, the problem of computing derivatives is severely ill-posed. Even a minor amount of noise in the data can result in significant errors in computation, see e.g., \cite{Ahnert:cpc2007, Engl:Kluwer1996}.
We present a well-known example in 1D. Consider $\eta(x) = \sin(nx)$, where $n$ is an integer and $x\in \mathbb{R}$. 
We have 
\begin{equation}
	f^{\delta}(x) = f^*(x) + \delta \sin(nx), 
	\quad {f^{\delta}}'(x) = {f^*}'(x) + n\delta \cos(nx),
	\label{sin_ex}
\end{equation}
 for $x\in \mathbb{R}$. 
The absolute error in this example is $\|n\delta \cos(nx)\|_{L^{\infty}(\R)} = n \delta \gg 1$ as $n \to \infty$, even though the perturbation term $\delta \sin(nx)$ is compatible to the noise level $\delta$, which is assumed to be small.
Due to the importance of differentiating noisy data, much research has been conducted in this area. In the following, we review some commonly employed methods for carrying out this task. 
\begin{enumerate}
	\item 
The finite difference method is the most natural approach for differentiation. An approximation of the derivative of $f^{\delta}$ can be obtained using the following formula:
\begin{equation}
{f^{\delta}}'(x) \approx \frac{f^{\delta}(x + h) - f^{\delta}(x)}{h} 
\approx {f^*}'(x)+  \frac{\delta(\eta(x + h) -\eta(x))}{h}
\label{FD}
\end{equation}
for $x \in \R$ where $h$ is a step size.
The step size $h$, which serves as a regularization factor, must be chosen appropriately to ensure stability. It is important to note that if the step size $h$ is too small, the expression on the right-hand side of \eqref{FD} becomes too large, leading to significant computational errors.
This is because $\eta(x + h)$ and $\eta(x)$ might take any value and, hence, the term $\frac{\delta(\eta(x + h) -\eta(x))}{h}$ might diverge as $h \to 0$.
 For further details and the strategy to choose the step size $h$, we refer the reader to \cite{Groetsch:amm1991, Ramm:mc2001}.

\item Tikhonov optimization is a more commonly used approach that involves introducing cost functionals with a Tikhonov regularization term. An example of such a cost functional is
\begin{equation}
J(u) = \Big\|f^{\delta}(x) -f^{\delta}(-R) - \int_{\R}^x u(s)ds\Big\|^2_{L^2(-R, R)} + \epsilon\|u\|_{X}^2,
\label{knowles}
\end{equation}
where $\epsilon > 0$ is the regularization parameter and $X$ is a functional space.
The absolute minimum of this cost functional provides an approximation to ${f^*}'$. This approach is discussed in more details in \cite{Knowles:ejde2014}.  Selecting the appropriate regularization parameter $\epsilon$ for each functional space $X$ for the regularization term in \eqref{knowles} is not a trivial task. A well-known method for determining the regularization parameter is Morozov's discrepancy principle, which is detailed in \cite{Engl:Kluwer1996, Morozov:springer1984, Ramlau:jnfao2002, Scherzer:SIAMJna1993}.

\item
An alternative method for computing derivatives of data involves the use of cubic spline curves to smooth the data. It is worth mentioning that the cubic smoothing spline was initially introduced in \cite{Reinsch:nm2967, Reinsch:nm1971, Schoenberg:PAMS1964}. Suppose that $f^\delta$ is the given data, and let $S$ be the natural cubic spline that approximates $f^\delta$. Here, by natural cubic spline, we mean that the function $S$ is twice continuously differentiable with $S''(-R) = S''(R) = 0$, and $S$ coincides on each subinterval $[x_{i - 1}, x_i]$ with some cubic polynomial, where $x_0 = -R < x_1 < x_2 < \dots < x_n = R$, $n > 1$, is a partition of $(-R, R)$. 
One can then use the derivatives of $S$ to approximate the derivatives of $f^{*}$.
For a method that combines cubic spline and Tikhonov optimization technique to differentiate noisy data, we refer the reader to \cite{Hanke2001}.

\item Other methods for differentiation are discussed in \cite{Breugel:IEEEAccess, Friedrichs:tams1944, Knowles:ejde2014, Knowles:nm1995} and their references.

\end{enumerate}

The example of $f^{\delta}(x) = f^*(x) + \delta \sin(nx)$ in \eqref{sin_ex} highlights the ill-posed nature of the problem, which arises from the high-frequency components of the data. It suggests cutting off these components to increase stability.
More precisely, given a function $f^{\delta} \in L^2(\Omega)$, we approximate
\begin{equation}
	f^{\delta}(\x) = \sum_{n = 1}^\infty f_n^{\delta}\Phi_n(\x) \simeq \sum_{n = 1}^N f_n^{\delta}\Phi_n(\x)
	\quad
	\x \in \Omega
	\label{1.5}
\end{equation}
for some cut-off number $N$,
where $\{\Phi_n\}_{n = 1}^\infty$ is an orthonormal basis of $L^2(\Omega)$ and
\[
	f_n^{\delta}(\x) = \int_{\Omega} f^{\delta}(\x)\Phi_n(\x)d\x.
\]
Then, we can approximate the desired derivatives using the following formulas
\begin{equation}
	\nabla f^{\delta}(\x) = \sum_{n = 1}^N f_n^{\delta} \nabla \Phi_n(\x),
	\quad
	D^2 f^{\delta}(\x) = \sum_{n = 1}^N f_n^{\delta} D^2 \Phi_n(\x),
	\label{1.2}
\end{equation}
for $\x \in \Omega.$
Although the numerical method based on \eqref{1.2} may not be rigorously valid due to the term-by-term differentiation of the series in \eqref{1.5}, it is an effective technique for numerical computation. Proving the analytical validity of this approach as $N$ approaches $\infty$ poses a significant challenge, which is not addressed in our paper. Nevertheless, our method remains acceptable. In fact, it is well-known that Problem \ref{p} is highly unstable, and therefore, a regularized technique is necessary. Truncating the series in \eqref{1.5} serves as an option for this purpose.

The main concern in \eqref{1.2} is to identify a suitable basis $\{\Phi_n\}_{n \geq 1}$ that enables reliable numerical differentiation through \eqref{1.2}.
We observe that using a popular basis such as Legendre polynomials or trigonometric functions in the well-known Fourier expansion might not be appropriate. 
This is mainly due to the fact that the initial function of these bases is a constant, which has an identically zero derivative. Consequently, the contribution of $f_1$ in \eqref{1.2} is overlooked. The lack of this information significantly reduces the accuracy of the method.
The basis $\{\Phi_n\}_{n \geq 1}$ we use in this paper is constructed by extending the polynomial-exponential basis initially introduced by Klibanov, the fourth author of this paper, in \cite{Klibanov:jiip2017} and \cite[Section 6.2.3]{KlibanovLiBook}. 
The original goal of Klibanov was to use this basis for the numerical solution of some coefficient inverse problems, see, e.g. \cite{VoKlibanovNguyen:IP2020, Khoaelal:IPSE2021, KhoaKlibanovLoc:SIAMImaging2020, 
KlibanovLeNguyen:SIAM2020, KlibanovNguyen:ip2019, LeNguyen:JSC2022}, \cite[Chapters 7,10-12]{KlibanovLiBook}, and \cite{KlibanovTimonov:cac2023}.
This basis is an appropriate choice for \eqref{1.2} since it satisfies the necessary condition that derivatives of $\Phi_n$, $n \geq 1$, are not identically zero.
We will demonstrate the effectiveness of the polynomial-exponential basis by comparing some numerical results due to \eqref{1.2} in two cases when 
 $\{\Phi_n\}_{n \geq 1}$ is the polynomial-exponential basis 
and $\{\Phi_n\}_{n \geq 1}$ is the widely-used trigonometric basis.
 Additionally, we will compare the numerical results obtained by our method with those obtained by the Tikhonov regularization and the cubic spline methods, which will serve as further evidence of the efficiency of our approach.

The polynomial-exponential basis $\{\Phi_n\}_{n \geq 1}$ has an additional significant property, in addition to the "non-identically zero" characteristic mentioned earlier. This property states that for all values of $N \geq 2$, the matrix $(\langle \Phi_n, \Phi_m'\rangle_{L^2})_{1\leq m, n \geq N}$ is invertible. This property guarantees a global convergence of Carleman-based convexification methods to solve inverse scattering problems in \cite{VoKlibanovNguyen:IP2020, Khoaelal:IPSE2021, KhoaKlibanovLoc:SIAMImaging2020,KlibanovLiBook, LeNguyen:JSC2022}.

The organization of this paper is as follows. Section \ref{section_comp} provides an overview of our approach. Section \ref{sec_num} presents numerical results for one and two dimensions. A comparison between our method and other numerical approaches can be found in section \ref{sec5}. Section \ref{sec6} is for concluding remarks.

\section{The numerical method based on the truncation technique} \label{section_comp}

In \cite{Klibanov:jiip2017}, Klibanov proposed a special orthonormal basis of $L^2(-R, R)$ to numerically solve inverse problems. 
For any $x \in (-R, R)$, we define $\phi_n(x) = x^{n - 1}e^x$. It is clear that the set $\{\phi\}_{n \geq 1}$ is complete in $L^2(-R, R)$. 
By applying the Gram-Schmidt orthonormalization process to this set, we obtain an orthonormal basis for $L^2(-R, R)$, which is denoted by $\{\Psi_n\}_{n \geq 1}$ and named the polynomial-exponential basis. The basis $\{\Psi_n\}_{n \geq 1}$ plays a vital role in our method for addressing Problem \ref{p}. The product of $\{\Psi_n\}_{n \geq 1}$ is the basis we use to approximate the noisy data $f^{\delta}$.
For $\bn = (n_1, n_2, \dots, n_d)$, we define
\begin{equation}
	P_{\bn}(\x) = \Psi_{n_1}(x_1) \Psi_{n_2}(x_2) \dots \Psi_{n_d}(x_d)
	\quad 
	\mbox{for all }
	\x = (x_1, \dots, x_d) \in \Omega.
	\label{basis}
\end{equation}
 It is obvious that $\{P_{\bn}\}_{\bn \in \N^d}$ forms an orthonormal basis of $L^2(\Omega).$
For any function $f \in L^2(\Omega)$, we can write
\begin{equation}
	f(\x) = \sum_{\bn \in \N^d} a_{\bn} P_{\bn}(\x)
	= \sum_{\bn = (n_1, \dots, n_d) \in \N^d} a_{\bn} \Psi_{n_1}(x_1) \Psi_{n_2}(x_2) \dots \Psi_{n_d}(x_d)
	\quad 
	\x \in \Omega,
	\label{2.2}
\end{equation} 
where
\begin{equation}
	a_{\bn} = a_{n_1 n_2 \dots n_d} = \int_{\Omega} f(\x) P_{\bn}(\x)d\x
	= \int_{\Omega} f(x_1, x_2, \dots, x_d)  \Psi_{n_1}(x_1) \Psi_{n_2}(x_2) \dots \Psi_{n_d}(x_d)(\x)d\x.
\end{equation}
Motivated by the Galerkin approximation, we truncate the representation of the function $f$ in  \eqref{2.2} as
\begin{equation*}
	f(\x) 
	= \sum_{n_1 = 1}^\infty \dots \sum_{n_d = 1}^{\infty} 
	a_{\bn} \Psi_{n_1}(x_1)
	\Psi_{n_2}(x_2)
	\dots
	\Psi_{n_d}(x_d)
	\simeq
	\sum_{n_1 = 1}^{N_1} \dots \sum_{n_d = 1}^{N_d} 
	a_{\bn} \Psi_{n_1}(x_1)
	\Psi_{n_2}(x_2)
	\dots
	\Psi_{n_d}(x_d)
\end{equation*} 
for some positive integers $N_1, \dots, N_d$ where $\bn = (n_1, \dots, n_d)$ and $\x = (x_1, \dots, x_d)$.
Hence, the function $f$ approximately becomes
\begin{equation}
	f(\x) \simeq \sum_{n_1 = 1}^{N_1} \dots \sum_{n_d = 1}^{N_d} 
	a_{\bn} \Psi_{n_1}(x_1)
	\Psi_{n_2}(x_2)
	\dots
	\Psi_{n_d}(x_d).
	\label{2.4}
\end{equation}
\begin{remark}
	Truncating the Fourier series, as depicted in equation \eqref{2.4}, provides an effective regularization technique for addressing ill-posed problems like Problem \ref{p}. This approach involves eliminating the high-frequency components of function $f$, resulting in a partial reduction of the ill-posedness of the problem. This technique has been successfully employed in solving inverse problems, where it enables the computation of numerical solutions even when significant noise levels are present, as evidenced by \cite{VoKlibanovNguyen:IP2020, KlibanovLeNguyen:SIAM2020, KlibanovNguyen:ip2019}.
	In
particular, highly noisy experimental data are treated by this technique in 
\cite{VoKlibanovNguyen:IP2020} and \cite[Chapter 10]{KlibanovLiBook}.
	 The cut-off numbers $N_1, \dots, N_d$ are utilized as regularization parameters in our method, and their determination can be achieved through numerical methods using the available data without requiring any knowledge of the noise level. We refer to Section \ref{sec3.1} for further details on selecting these cut-off numbers.
\end{remark}
After smoothing the data using \eqref{2.4}, we can compute the first derivatives of $f$ via the formulas
\begin{equation}
	\partial_{x_i} f(\x) \simeq
	\sum_{n_1 = 1}^{N_1} \dots \sum_{n_d = 1}^{N_d} 
	a_{\bn} \Psi_{n_1}(x_1)
	\Psi_{n_2}(x_2)
	\dots
	\Psi_{n_i}'(x_i)
	\dots
	\Psi_{n_d}(x_d).
	\label{2.5}
\end{equation}
for $i \in \{1, \dots, d\}$.
The second derivatives of $f$ are computed using
\begin{equation}
	\partial_{x_i x_j} f(\x) 
	\simeq
	\left\{
	\begin{array}{ll}
	\ds\sum_{n_1 = 1}^{N_1} \dots \sum_{n_d = 1}^{N_d} 
	a_{\bn} \Psi_{n_1}(x_1)
	\Psi_{n_2}(x_2)
	\dots
	\Psi_{n_i}'(x_i)
	\dots
	\Psi_{n_j}'(x_j)
	\dots
	\Psi_{n_d}(x_d) & i \not = j,
	\\
	\ds\sum_{n_1 = 1}^{N_1} \dots \sum_{n_d = 1}^{N_d} 
	a_{\bn} \Psi_{n_1}(x_1)
	\Psi_{n_2}(x_2)
	\dots
	\Psi_{n_i}''(x_i)
	\dots
	\Psi_{n_d}(x_d)
	&i = j,
	\end{array}
	\right.
	\label{2.6}
\end{equation}
for $i, j \in \{1, \dots, d\}.$
The approximations in \eqref{2.5}  and \eqref{2.6} suggest Algorithm \ref{alg} to solve Problem \ref{p}.
We refer the reader to Section \ref{sec3.1} for a numerical method to choose the cut-off numbers in the first step of Algorithm \ref{alg}.

\begin{algorithm}[h!]
\caption{\label{alg}The procedure to compute derivatives of noisy data}
	\begin{algorithmic}[1]
	\State \label{s1} Choose cut-off numbers $N_1, N_2, \dots, N_d$.
	 \State \label{s2} Construct the orthonormal basis $\{\Psi_n\}_{n \geq 1}$ defined in the first paragraph of Section \ref{section_comp}.
	 \State \label{s3}
	 For each $\bn \in \{(n_1, \dots, n_d): 1 \leq n_1 \leq N_1, \dots, 1 \leq n_d < N_d\}$,	compute the coefficients $a_{\bn}$
	 \begin{equation}
	 	a_{\bn} = a_{n_1 n_2 \dots n_d}  = \int_{\Omega} f(x_1, x_2, \dots, x_d) \Psi_{n_1}(x_1) \Psi_{n_2}(x_2) \dots \Psi_{n_d}(x_d) d\x.
	 \end{equation}
	\State \label{s4} Find the first derivatives of $f$ via formulas \eqref{2.5}.
	\State \label{s24} Find the second derivatives of $f$ via formulas \eqref{2.6}.
\end{algorithmic}
\end{algorithm}

\begin{remark}
	The selection of the polynomial-exponential basis $\{P_{\bf n}\}_{{\bf n} \in \N^d}$ is critical for achieving success in our approach. 
	Upon revisiting equations \eqref{2.5} and \eqref{2.6}, it becomes apparent that the accuracy of the computed results is heavily reliant on the contribution of the Fourier coefficients $a_{\bf n}$, ${\bf n} \in \N^d$. 
	Therefore, if an inappropriate basis such as $\{Q_{\bf n}\}_{n \geq 1}$ is chosen, which contains an entry $Q_{\bf m}$ with zero derivative for some ${\bf m} \in \N^d$, then the corresponding Fourier coefficient $b_{\bf m} = \int_{\Omega} f(\x) Q_{\bf m}(\x) d\x$ would not be present in equations \eqref{2.5} and \eqref{2.6}. 	
	For instance, in 1D, the first entry of the trigonometric basis involving $\cos\Big(\frac{\pi n x}{R}\Big)$, $n \geq 0$, and $\sin\Big(\frac{\pi m x}{R}\Big)$, $m \geq 1$, used in the Fourier expansion is the constant $\cos 0 = 1$,
	whose derivative vanishes.
	 Consequently, the first term in equations \eqref{2.5} and \eqref{2.6} for computed derivatives will be missing.
	 This leads to a decrease in accuracy during computation as the first term in the series could be the most crucial one.
	 Thus, the trigonometric basis is unsuitable for our method. This observation will be confirmed through numerical experiments in this paper.
\end{remark}

\section{Numerical examples} \label{sec_num}

We implement Algorithm \ref{alg} for the cases of 1D and 2D. 
In all tests below, we compute the derivatives of a function $f^*$ using the noisy data given by
\begin{equation}
f^{\delta}(\x) = f^*(\x)(1 + \delta \mbox{rand}(\x))
\end{equation}
where rand$(\x)$ is a function that generates uniformly distributed random numbers in $[-1, 1].$ The noise levels are set to be $5\%$, $10\%,$ and $20\%.$ 

\subsection{Choosing the cut-off number} \label{sec3.1}
\begin{figure}[ht]
	\subfloat[The data $f^{\delta}$ (dash-dot) and its approximation $\sum_{n = 1}^N f_n^\delta \Psi_n(x)$ (solid) when $N = 10$.]{\includegraphics[width=.3\textwidth]{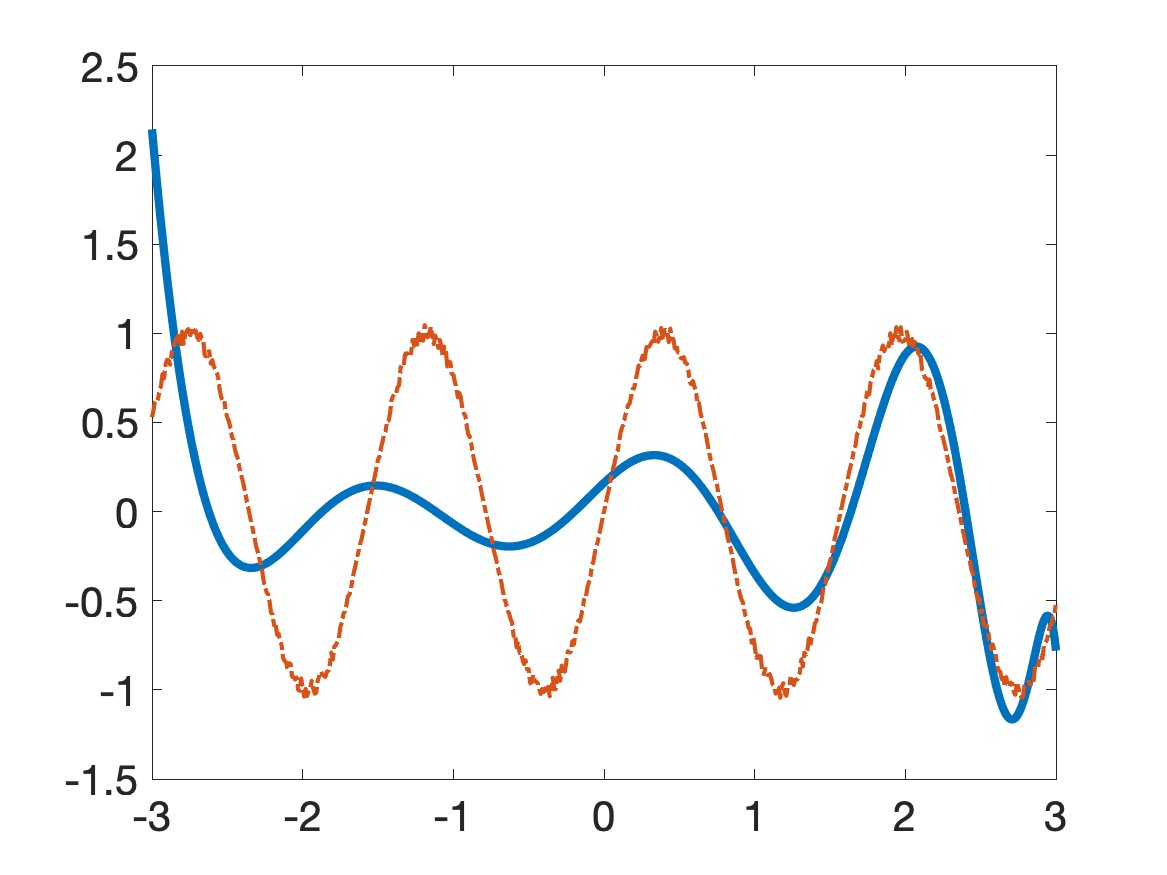}}
	\quad
	\subfloat[The data $f^{\delta}$ (dash-dot) and its approximation $\sum_{n = 1}^N f_n^\delta \Psi_n(x)$ (solid) when $N = 15$.]{\includegraphics[width=.3\textwidth]{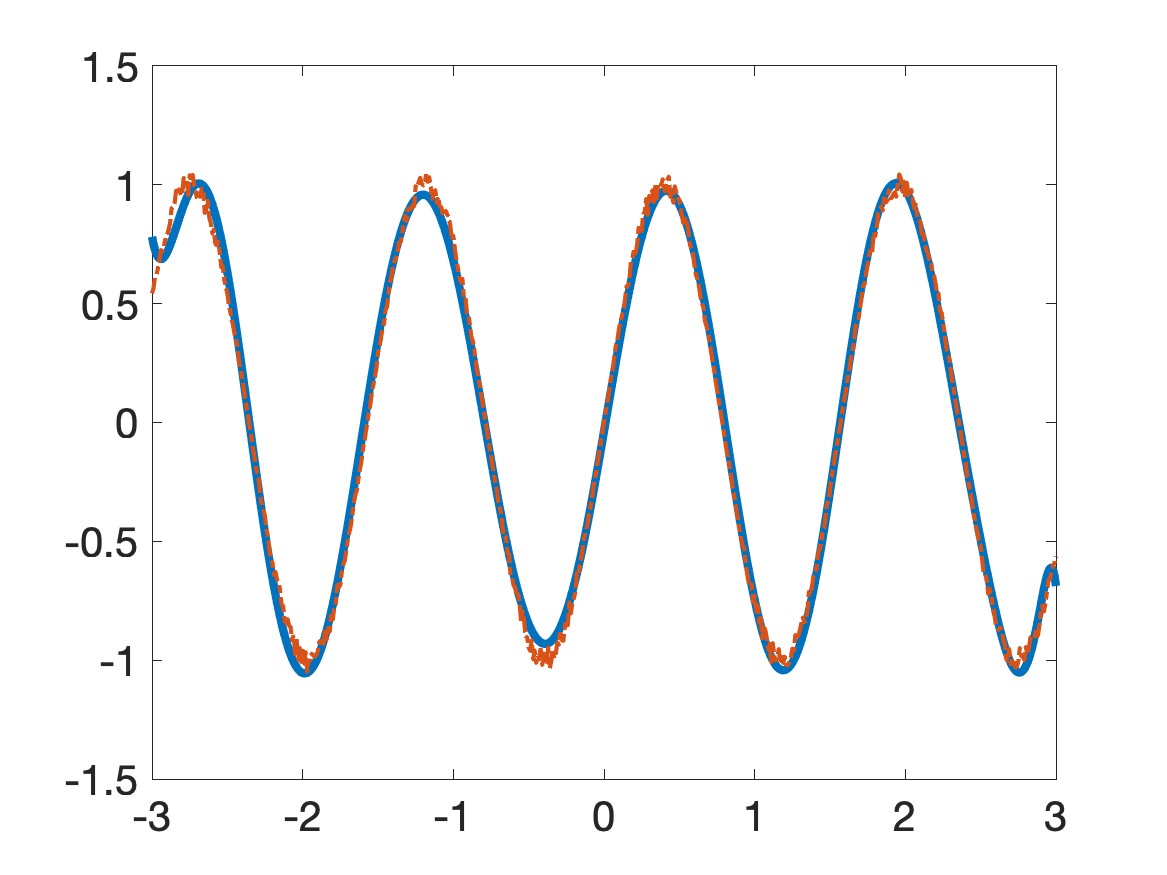}}
	\quad
	\subfloat[The data $f^{\delta}$ (dash-dot) and its approximation $\sum_{n = 1}^N f_n^\delta \Psi_n(x)$ (solid) when $N = 20$.]{\includegraphics[width=.3\textwidth]{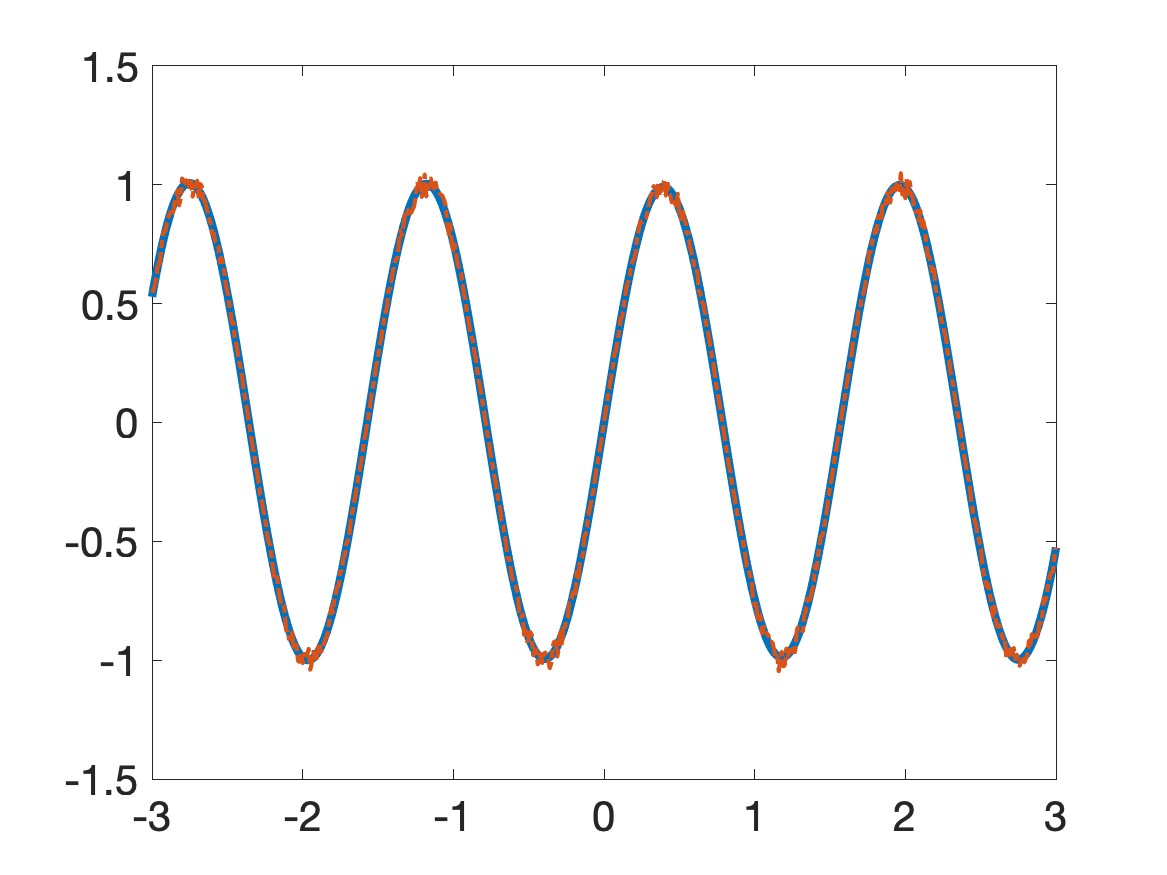}}
	
	\subfloat[The function $\varphi_N$ when $N = 10$.]{\includegraphics[width=.3\textwidth]{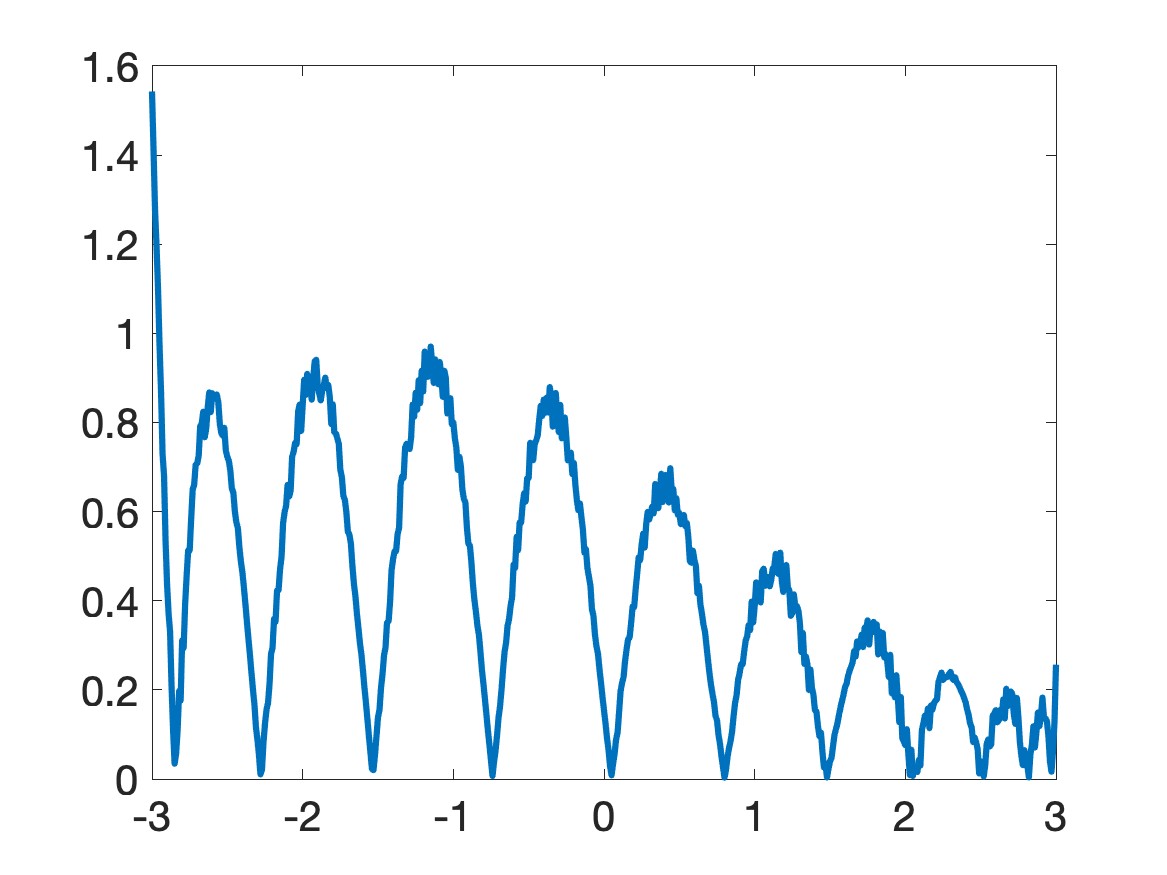}}
	\quad
	\subfloat[The function $\varphi_N$ when $N = 15$.]{\includegraphics[width=.3\textwidth]{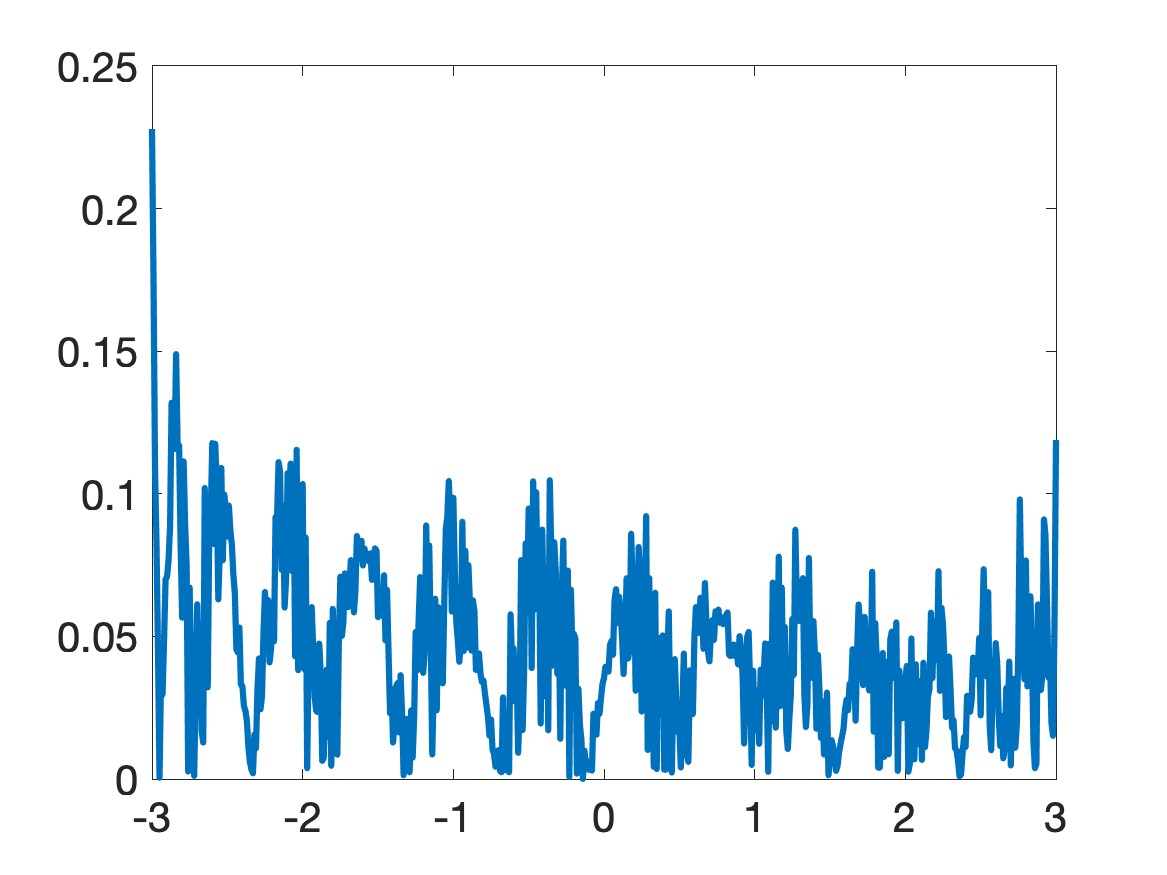}}
	\quad
	\subfloat[The function $\varphi_N$ when $N = 20$.]{\includegraphics[width=.3\textwidth]{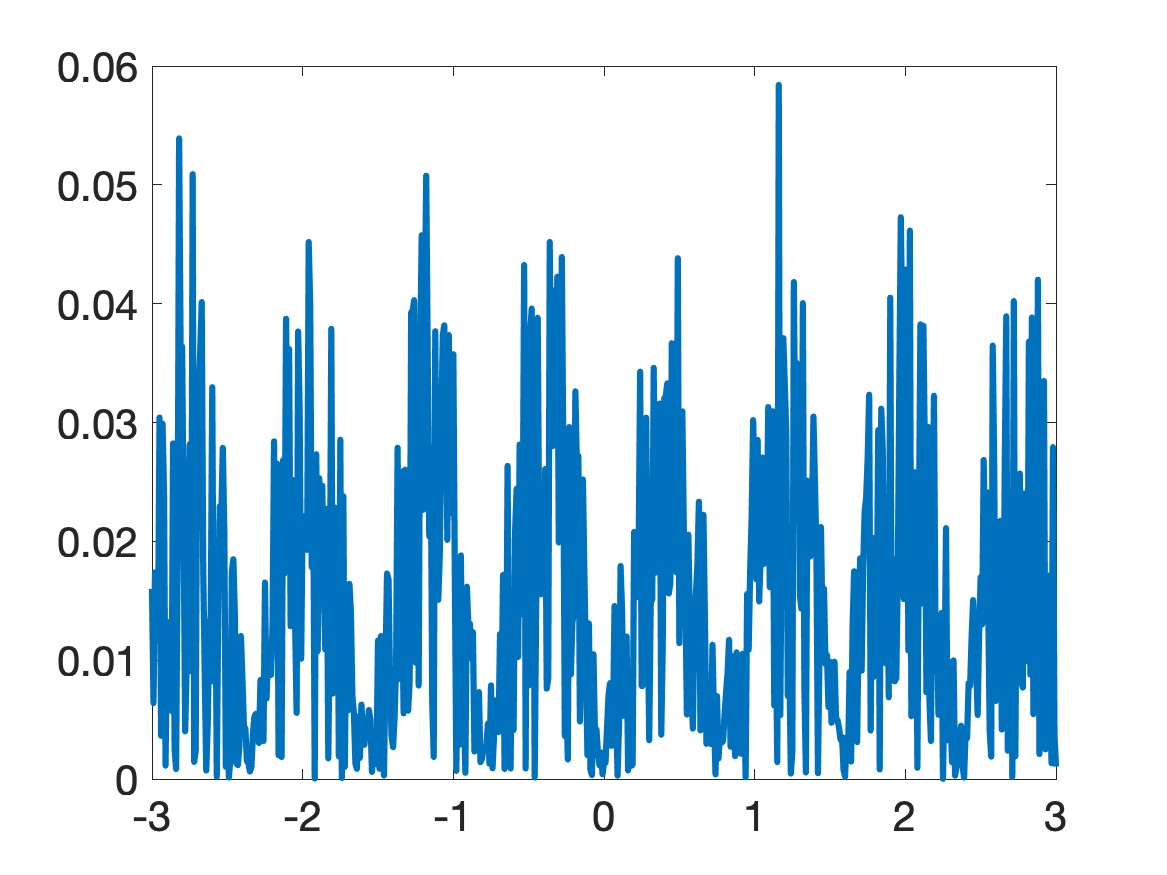}}
	\caption{\label{chooseN} An example of how to choose the cut-off number $N$ where the data, with $5\%$ noise, is taken from Test 1 below. It is evident that $N = 20$ should be chosen for Test 1.}
\end{figure}

Before presenting the numerical results, we explain the process of selecting the cut-off number in Step \ref{s1} of Algorithm \ref{alg}. To simplify matters, we will only present our strategy for the one-dimensional case. The same approach is utilized for selecting the cut-off numbers in the two-dimensional case.
For each $N > 1$, we compare noisy data $f^{\delta}$ with the sum $ \sum_{n = 1}^N f_n^{\delta} \Psi_n(x)$. In other words, for each given data $f^{\delta}$, we compute 
\[
	\varphi_N(x) = \frac{\big|f^{\delta}(x) - \sum_{n = 1}^N f_n^\delta \Psi_n(x)\big|}{\|f^{\delta}\|_{L^{\infty}}},
	\quad 
	\mbox{where } 
	f_n^{\delta} = \int_{-R}^R f^{\delta}(x)\Psi_n(x)dx,
\]
for several values of $N$. Then we choose $N$ is $\|\varphi_N\|_{L^{\infty}(-R, R)}$ is sufficiently small.
An illustration of this strategy is provided in Figure \ref{chooseN}.

This methodology is implemented in all of the following tests. The same procedure is used for the experiment involving two dimensions.

\subsection{The case of one dimension}\label{sec4.1}
We compute the first and second derivatives of some functions $f^*(x)$ defined in $(-R, R)$ where $R = 3.$ 
On $(-R, R)$, we arrange a partition 
\begin{equation}
	\mathcal{G} = \big\{x_i = -R + (i-1)h, i = 1, \dots, N\big\} \subset [-3, 3]
	\label{3.2222}
\end{equation}
where the step size $h = 0.001$ and $N = 6001$.

\subsubsection{Test 1} 
In this test, we take $f^*(x) = \sin(4x)$ as our function of interest. The accurate first and second-order derivatives are $f_{\rm true}'(x) = 4\cos(x)$ and $f_{\rm true}''(x) = -16 \sin(x)$, respectively. Our chosen cut-off number for this test is $N = 20$. 
\begin{figure}[ht]
    \centering
    \subfloat[The true and computed $f'(x)$; noise level $\delta = 5\%$]{\includegraphics[width = .3\textwidth]{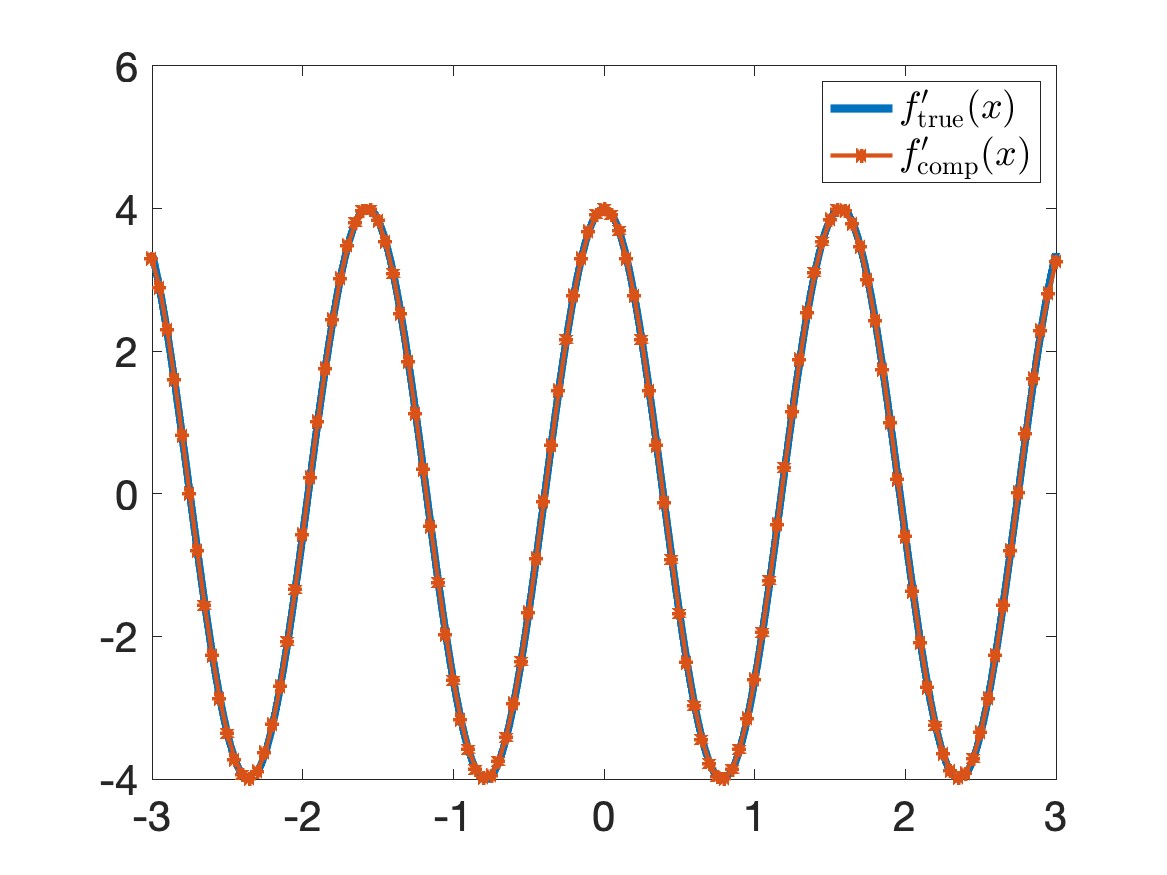}}
    \quad
    \subfloat[The true and computed $f'(x)$; noise level $\delta = 10\%$]{\includegraphics[width = .3\textwidth]{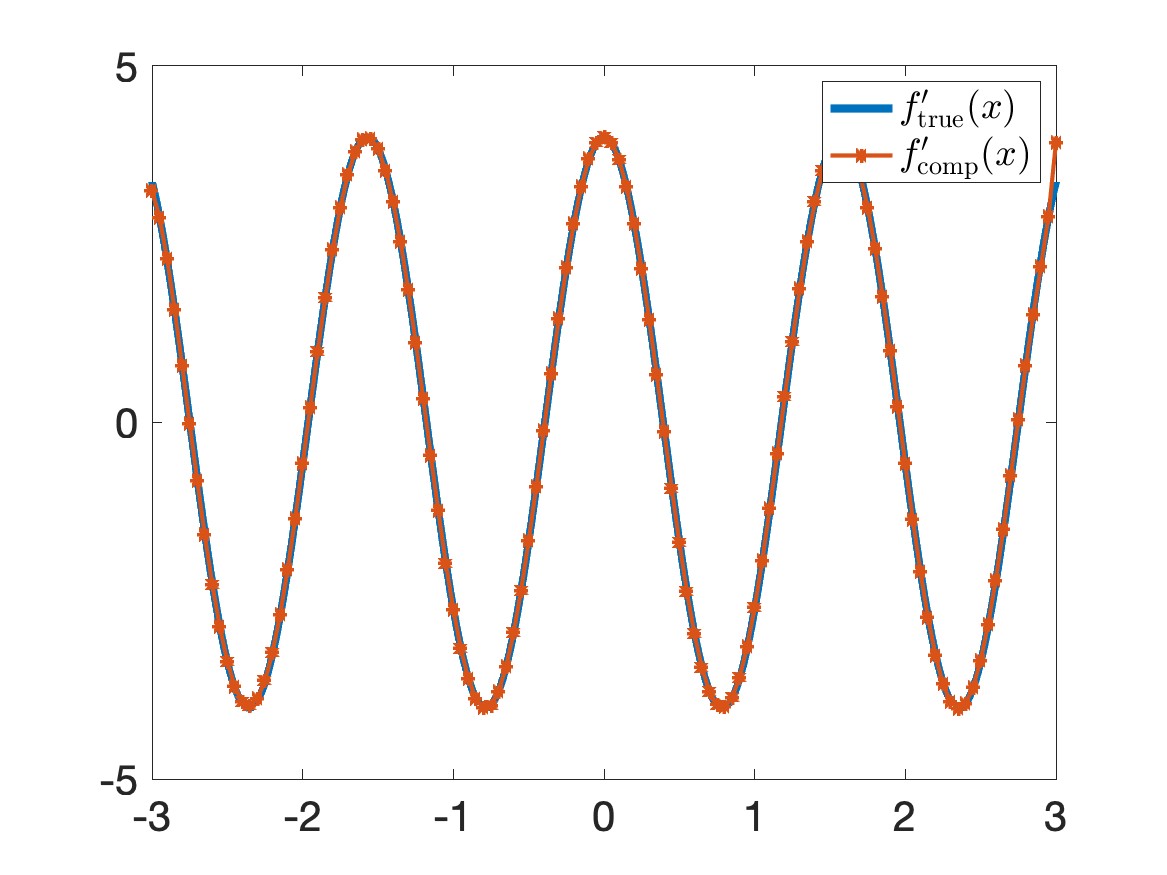}}
    \quad
    \subfloat[The true and computed $f'(x)$; noise level $\delta = 20\%$]{\includegraphics[width = .3\textwidth]{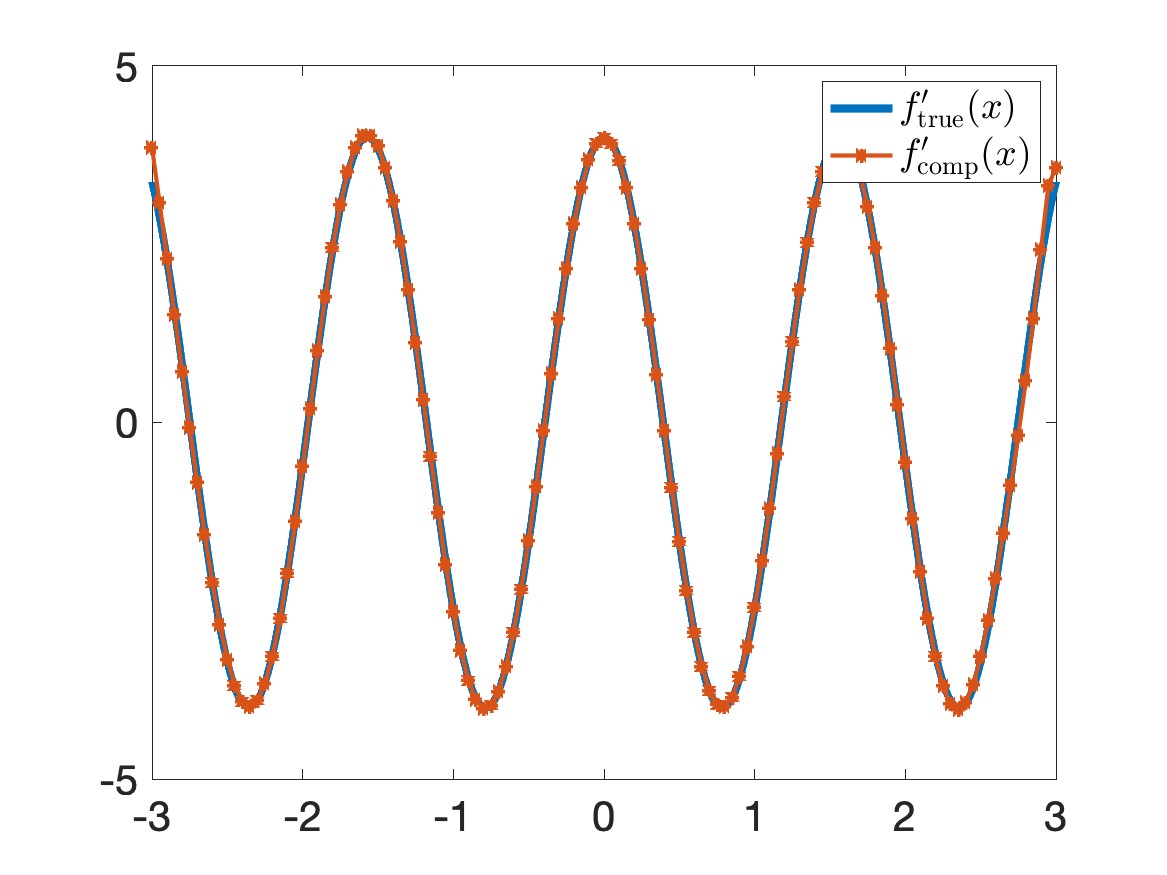}}

\subfloat[The true and computed $f''(x)$; noise level $\delta = 5\%$]{\includegraphics[width = .3\textwidth]{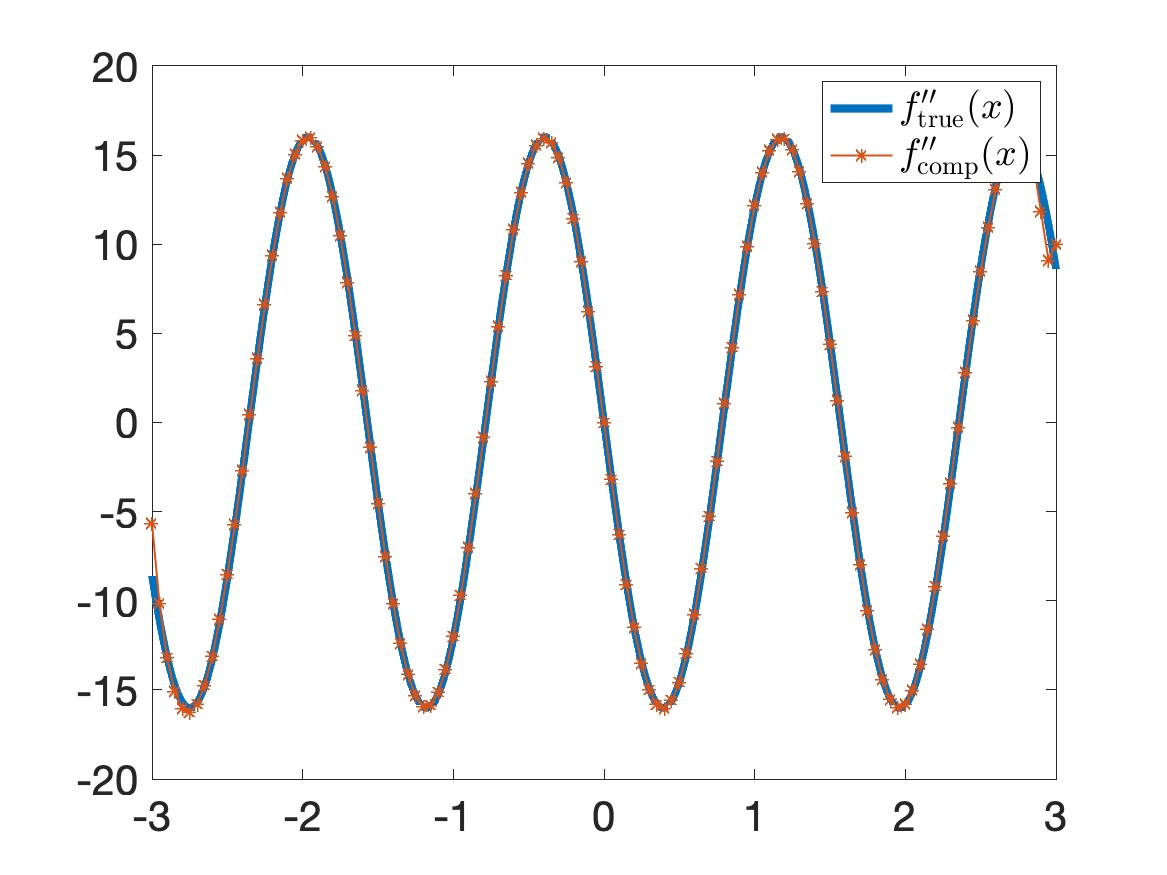}}
    \quad
    \subfloat[The true and computed $f''(x)$; noise level $\delta = 10\%$]{\includegraphics[width = .3\textwidth]{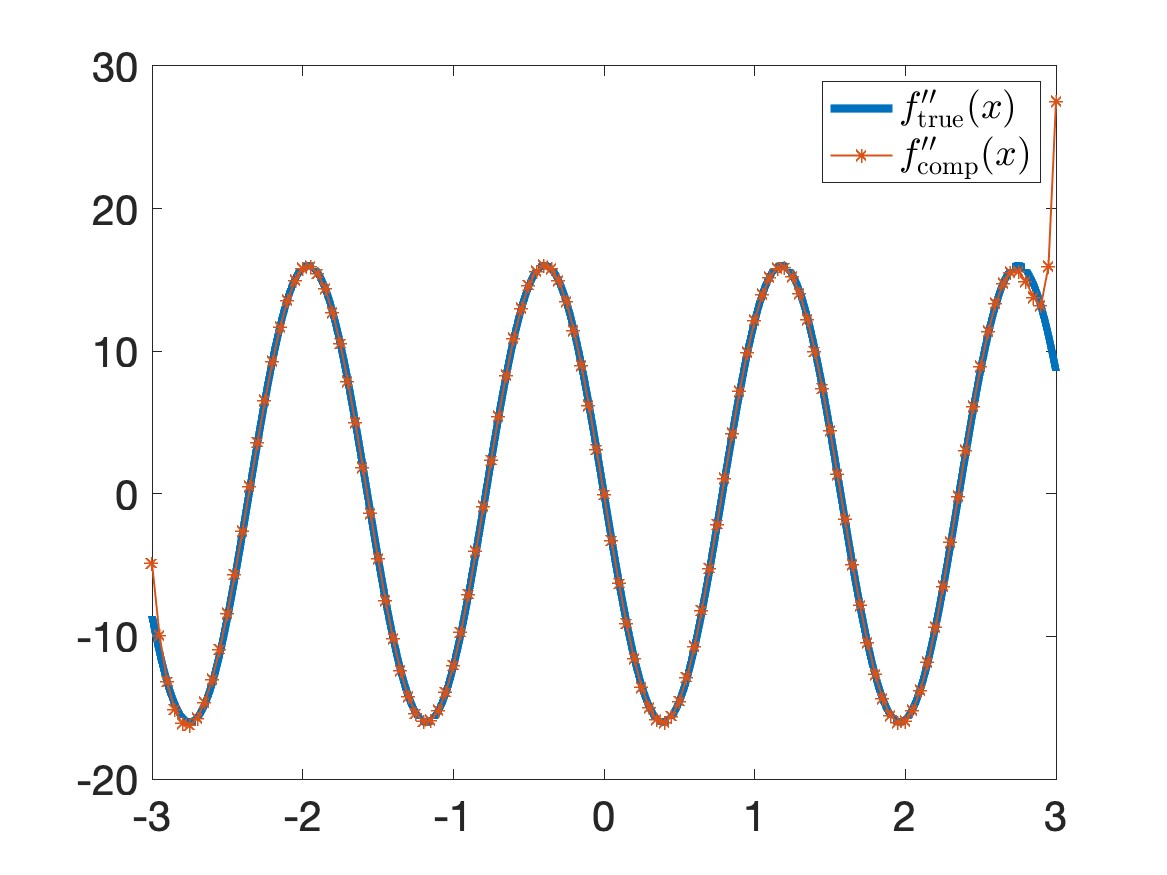}}
    \quad
    \subfloat[The true and computed $f''(x)$; noise level $\delta = 20\%$]{\includegraphics[width = .3\textwidth]{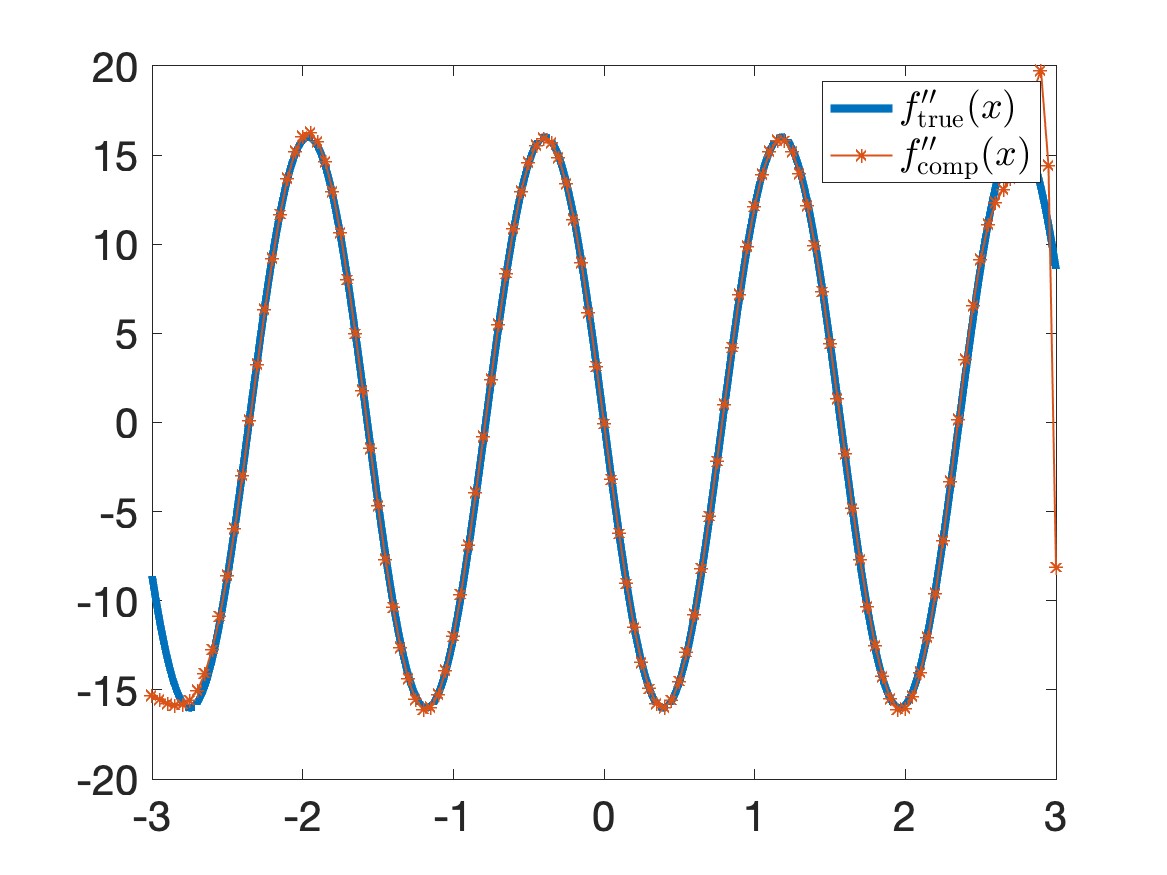}}
    
    \caption{\label{fig1} Test 1. The true and computed derivatives of the function $f(x) = \sin(4x)$. As can be seen, our method is highly proficient in accurately computing derivatives of both the first and second order. Notably, the errors are primarily concentrated at the endpoints of the computed domain.}   
\end{figure}
 
The graphical representation in Figure \ref{fig1} illustrates both the computed and true derivatives of $f^*$. Based on the outcomes shown in Figure \ref{fig1}, it is apparent that Algorithm \ref{alg} performs effectively in computing derivatives, even when there is a high level of noise, up to 20\%. It is worth mentioning that the errors mostly occur at the endpoints of the computational interval, while the computed function's accuracy within the interval exceeds expectations. The first two columns of Table \ref{tab1.1} display the computed errors over the entire interval $(-3,3)$, while the last two columns show the computed errors over the interval $(-2,2)$.
\begin{table}[ht]
    \centering
\begin{tabular}{|c|c|c|c|c|}
    \hline 
    		&$\frac{\|f'_{\rm true} - f'_{\rm comp}\|_{L^2(-3, 3)}}{\|f'_{\rm true}\|_{L^2(-3, 3)}}$
    		& $\frac{\|f''_{\rm true} - f''_{\rm comp}\|_{L^2(-3, 3)}}{\|f''_{\rm true}\|_{L^2(-3, 3)}}$
		&$\frac{\|f'_{\rm true} - f'_{\rm comp}\|_{L^2(-2, 2)}}{\|f'_{\rm true}\|_{L^2(-2, 2)}}$
		&  $\frac{\|f''_{\rm true} - f''_{\rm comp}\|_{L^2(-2, 2)}}{\|f''_{\rm true}\|_{L^2(-2, 2)}}$
  \\  \hline
    		$\delta = 0.05$
		&	 0.0060	
		& 	0.0268
		& 0.0030
		&0.0195
  \\  \hline
    		$\delta = 0.10$
		&0.0110
		&0.0996
		&0.0031
		&0.0201
	  \\  \hline
    		$\delta = 0.20$
		&0.0260
		&0.1123
		&0.0073
		&0.0282
	\\\hline
\end{tabular}
\caption{\label{tab1.1}Test 1. The relative computed errors with respect to the $L^2$ norms over the intervals $(-3, 3)$ and $(-2, 2)$. As can be observed, the relative errors are lower than the noise level.}
\end{table}

\subsubsection{Test 2}

\begin{figure}[ht]
    \centering
    \subfloat[The true and computed $f'(x)$; noise level $\delta = 5\%$]{\includegraphics[width = .3\textwidth]{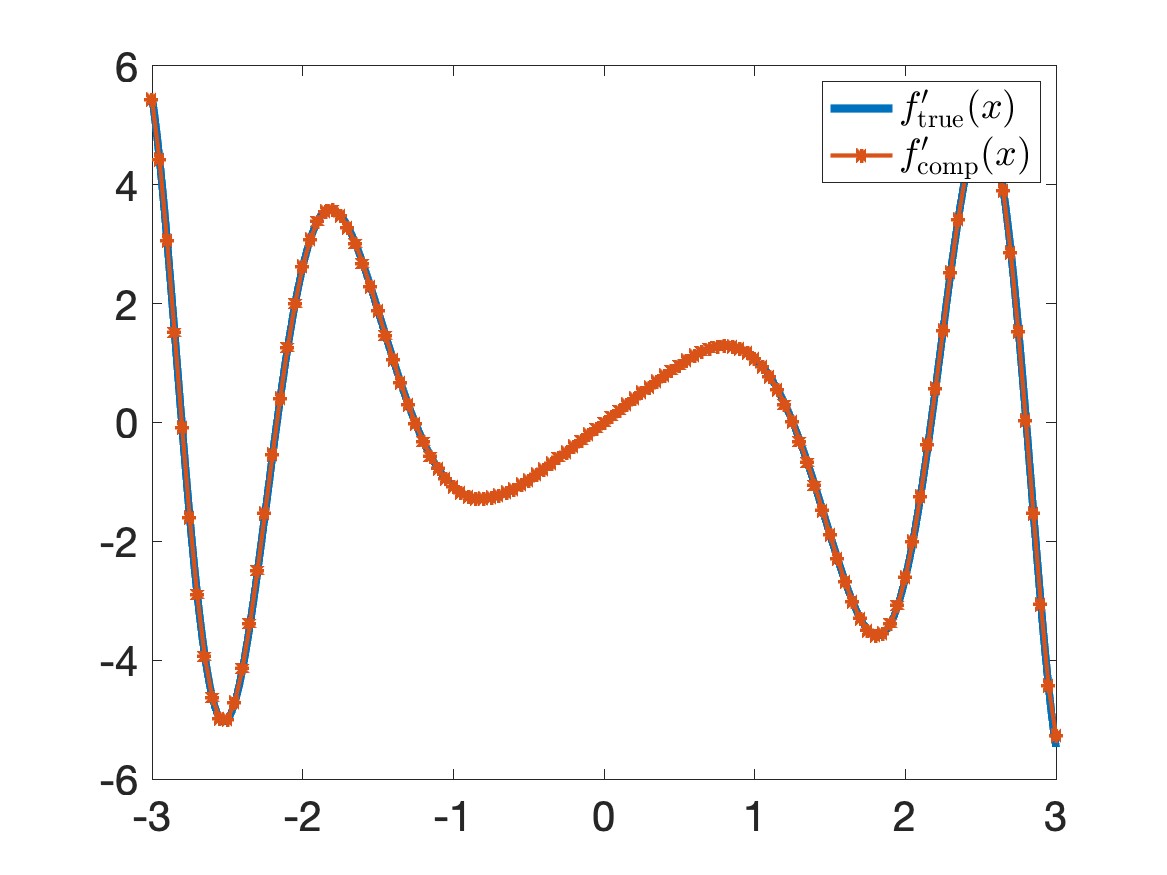}}
    \quad
    \subfloat[The true and computed $f'(x)$; noise level $\delta = 10\%$]{\includegraphics[width = .3\textwidth]{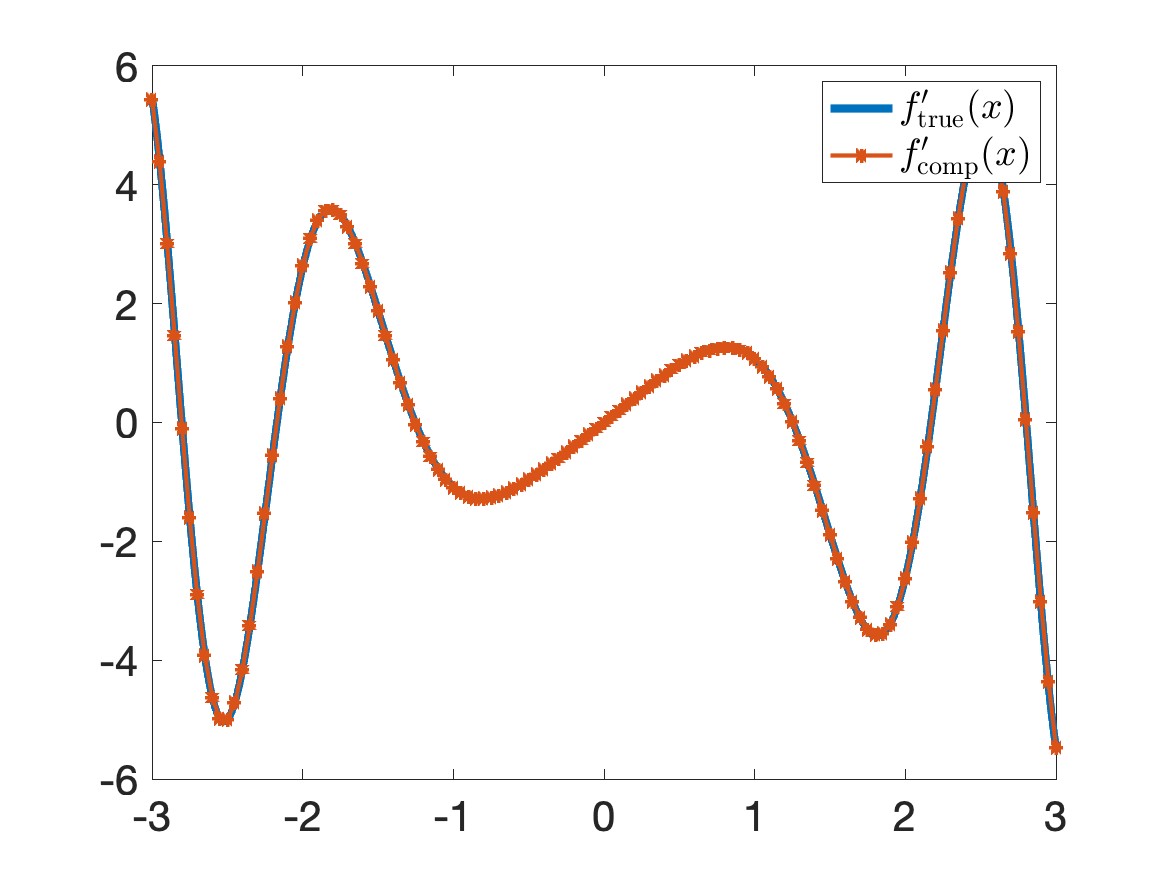}}
    \quad
    \subfloat[The true and computed $f'(x)$; noise level $\delta = 20\%$]{\includegraphics[width = .3\textwidth]{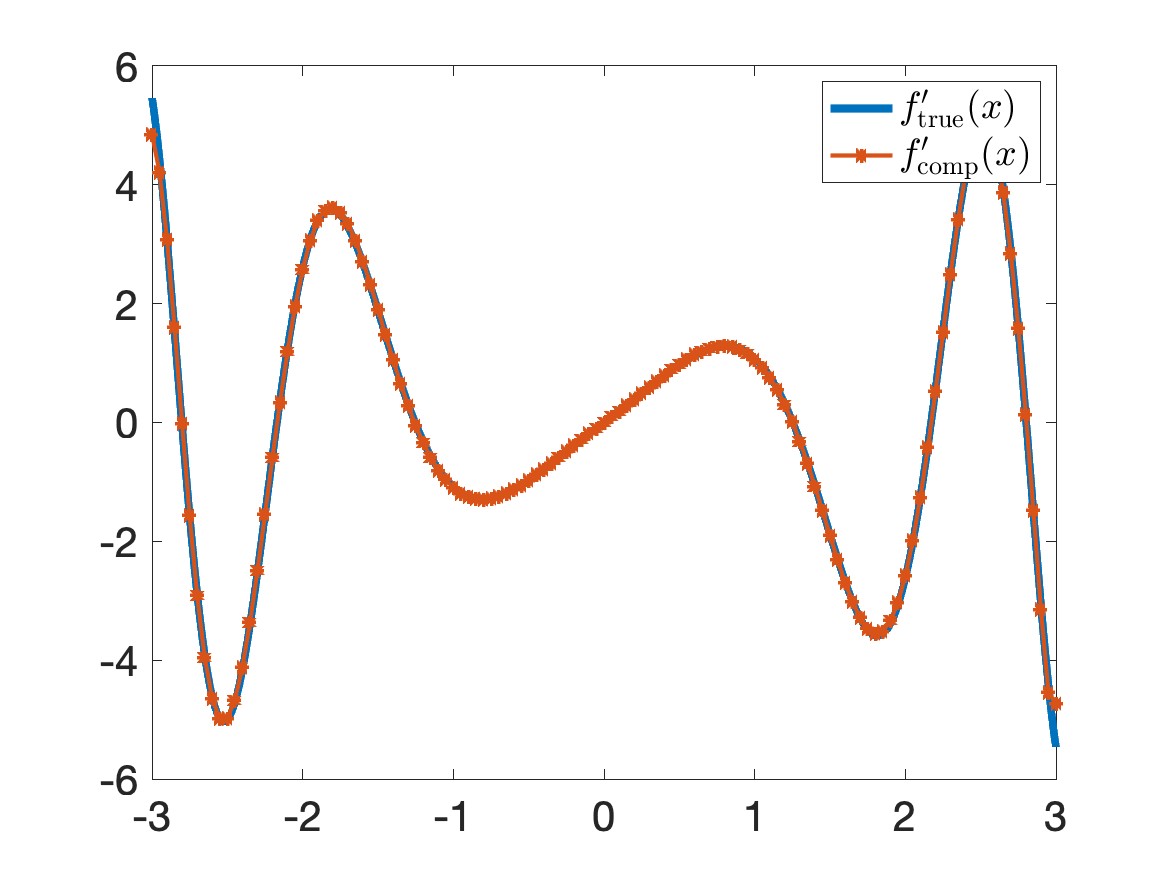}}

\subfloat[The true and computed $f''(x)$; noise level $\delta = 5\%$]{\includegraphics[width = .3\textwidth]{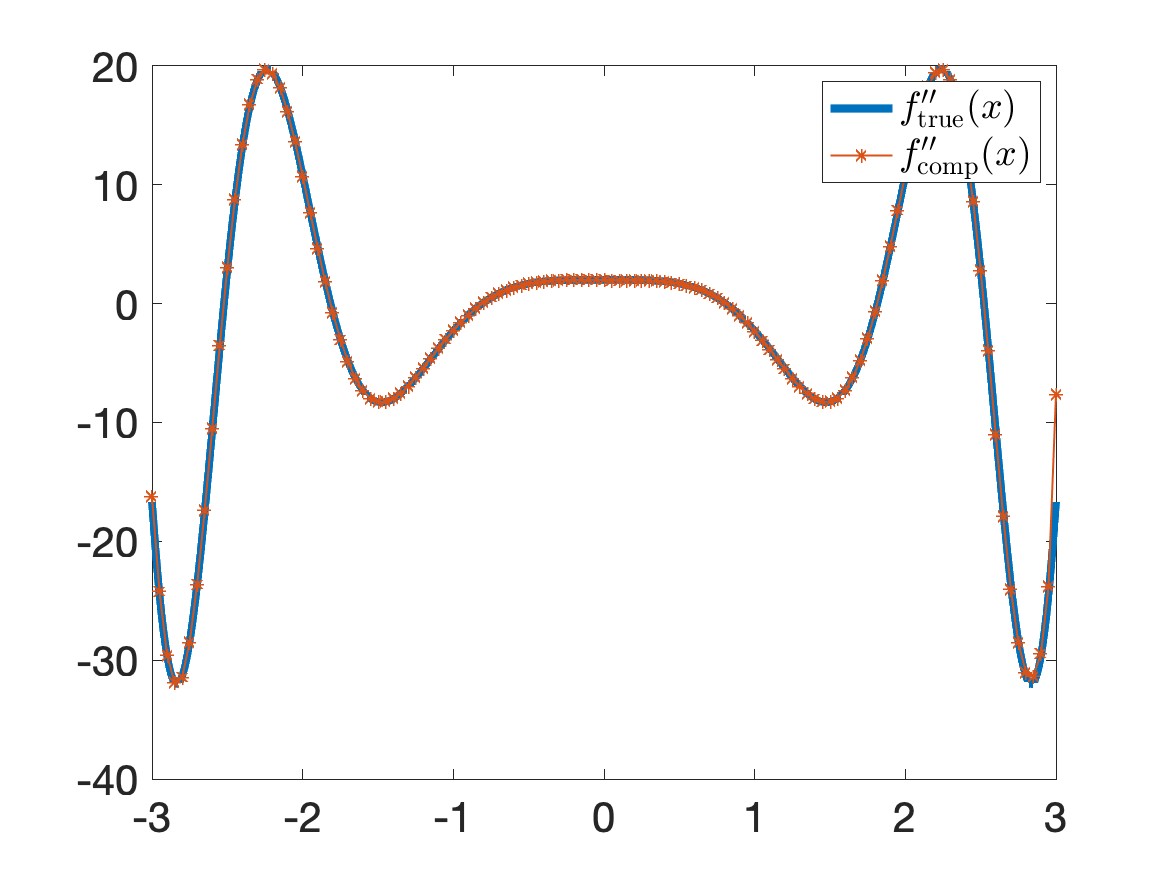}}
    \quad
    \subfloat[The true and computed $f''(x)$; noise level $\delta = 10\%$]{\includegraphics[width = .3\textwidth]{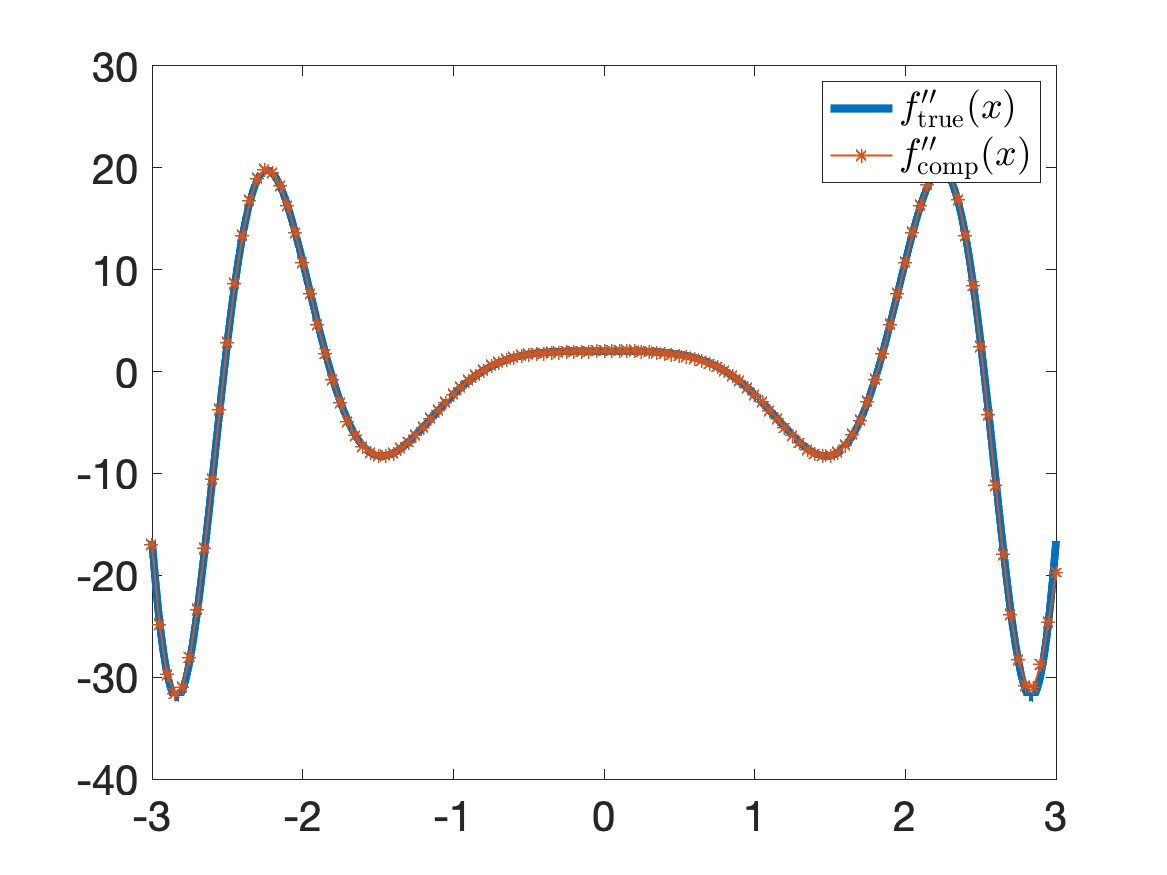}}
    \quad
    \subfloat[The true and computed $f''(x)$; noise level $\delta = 20\%$]{\includegraphics[width = .3\textwidth]{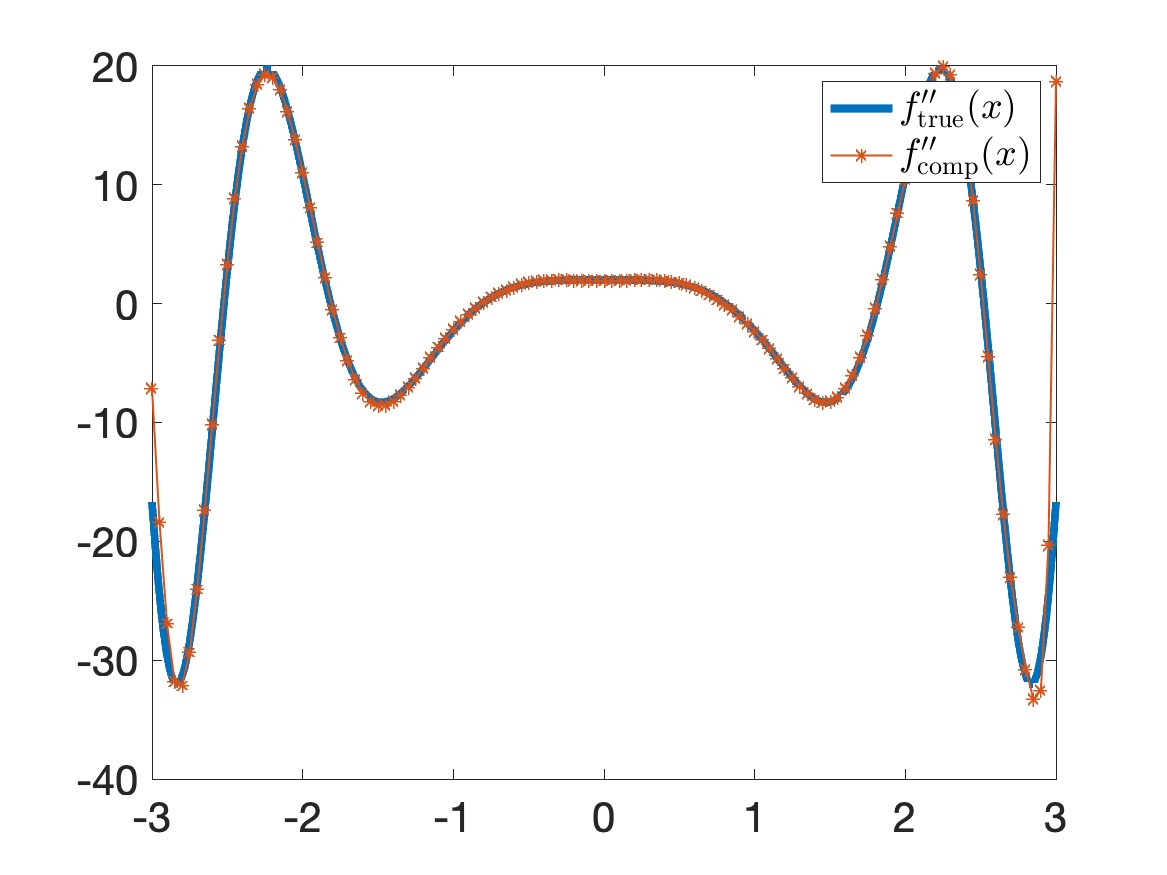}}
    
    \caption{\label{fig2} Test 2. The true and computed derivatives of the function $f(x) = \sin(x^2)$. As seen, the performance of Algorithm \ref{alg} in this test is out of expectation. The errors occur at the endpoints $\{-3, 3\}$ while the computed derivatives are accurate inside the interval $(-3,3)$.}   
\end{figure}

We calculate the derivatives for the function $f^*(x) = \sin(x^2)$. The true first-order and second-order derivatives for this function are given by $f_{\rm true}' = 2x \cos(x^2)$ and $f''_{\rm true} = 2 \cos(x^2) - 4x^2\sin(x^2)$, respectively. For this test, we use a cut-off number of $N = 25$. Figure \ref{fig2} displays the graphs of the computed derivatives.

Similar to Test 1, Figure \ref{fig2} visually confirms that our algorithms for computing the first and second derivatives are highly effective. The errors are primarily localized at the endpoints of the interval, whereas the computed function demonstrates noteworthy accuracy within the subinterval $(-2, 2)$. Detailed information can be found in Table \ref{tab2}.

\begin{table}[ht]
    \centering
\begin{tabular}{|c|c|c|c|c|}
    \hline 
    		&$\frac{\|f'_{\rm true} - f'_{\rm comp}\|_{L^2(-3, 3)}}{\|f'_{\rm true}\|_{L^2(-3, 3)}}$
    		& $\frac{\|f''_{\rm true} - f''_{\rm comp}\|_{L^2(-3, 3)}}{\|f''_{\rm true}\|_{L^2(-3, 3)}}$
		&$\frac{\|f'_{\rm true} - f'_{\rm comp}\|_{L^2(-2, 2)}}{\|f'_{\rm true}\|_{L^2(-2, 2)}}$
		&  $\frac{\|f''_{\rm true} - f''_{\rm comp}\|_{L^2(-2, 2)}}{\|f''_{\rm true}\|_{L^2(-2, 2)}}$
  \\  \hline
    		$\delta = 0.05$
		& 0.0052	 	
		& 0.0380	
		& 0.0017
		& 0.0309
  \\  \hline
    		$\delta = 0.10$
		&0.0074
		& 0.0955
		& 0.0047
		&0.0484
	  \\  \hline
    		$\delta = 0.20$
		&0.0240
		&0.1734
		&0.0117
		&0.0704
	\\\hline
\end{tabular}
\caption{\label{tab2} Test 2. The computed errors with respect to the $L^2$ norms over the intervals $(-3, 3)$ and $(-2, 2)$. It is apparent that the relative errors are below the noise level.}
\end{table}

\subsection{The case of two dimensions}

Our computation involves setting $\Omega = (-R, R)^2$, where $R = 3$. Similar to the one-dimensional case, we implement Algorithm \ref{alg} in the finite difference scheme. To define our grid, we set
\[
	\mathcal{G} = \big\{\x_{ij} = (x_i = -R + (i-1) h, y_j = -R + (j-1)h: 1\leq i, j\leq N_{\x}\big\}
\]
where $N_\x = 601$ and $h = 2R/(N_\x - 1)=0.01$. We choose the cut-off numbers $N_1$ and $N_2$ using the same strategy as in Section \ref{sec4.1}, i.e., we find $N_1$ and $N_2$ such that the approximation in \eqref{2.2} is satisfactory. For both tests presented below, we set $N_1 = N_2 = 20$.

\subsubsection{Test 3} 
\begin{figure}[ht]
	\subfloat[The true function $|\nabla f^*|$]{\includegraphics[width = .3\textwidth]{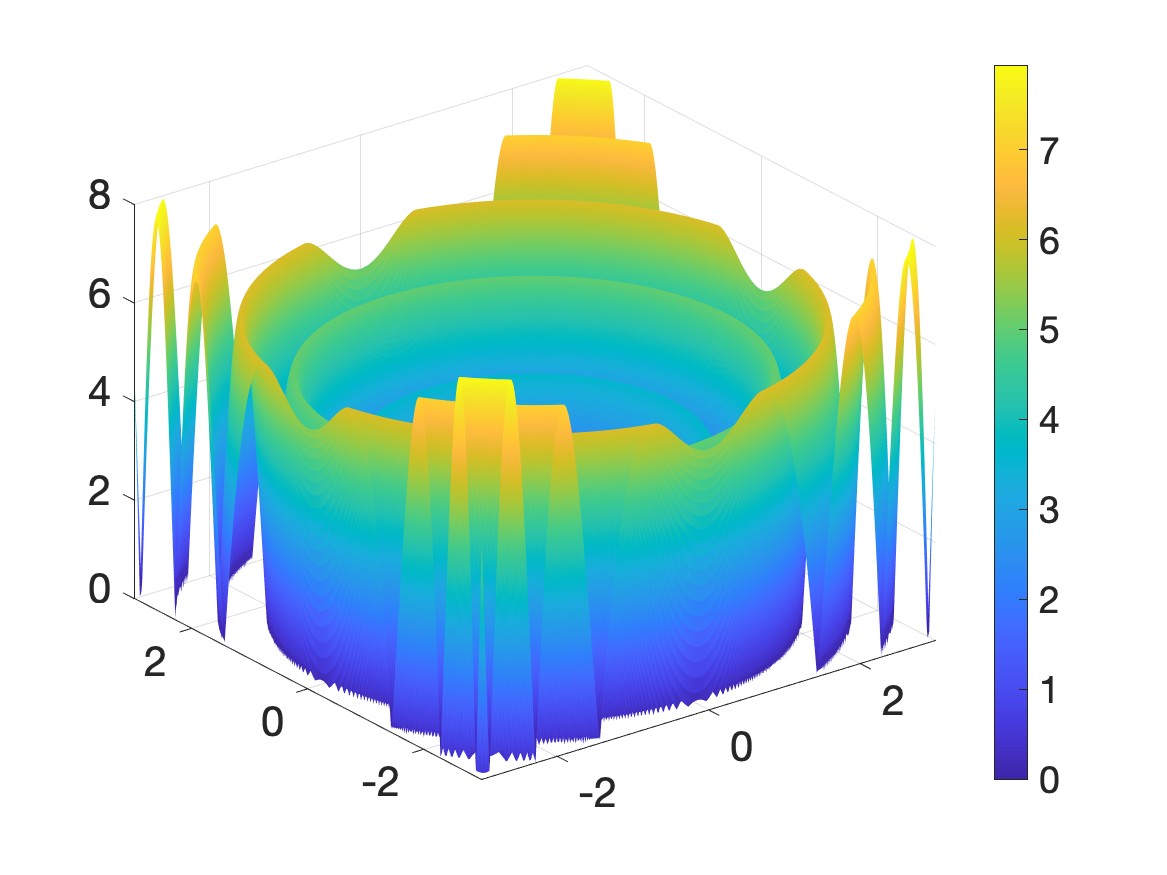}}
	\quad
	\subfloat[The computed function $|\nabla f_{\rm comp}|$]{\includegraphics[width = .3\textwidth]{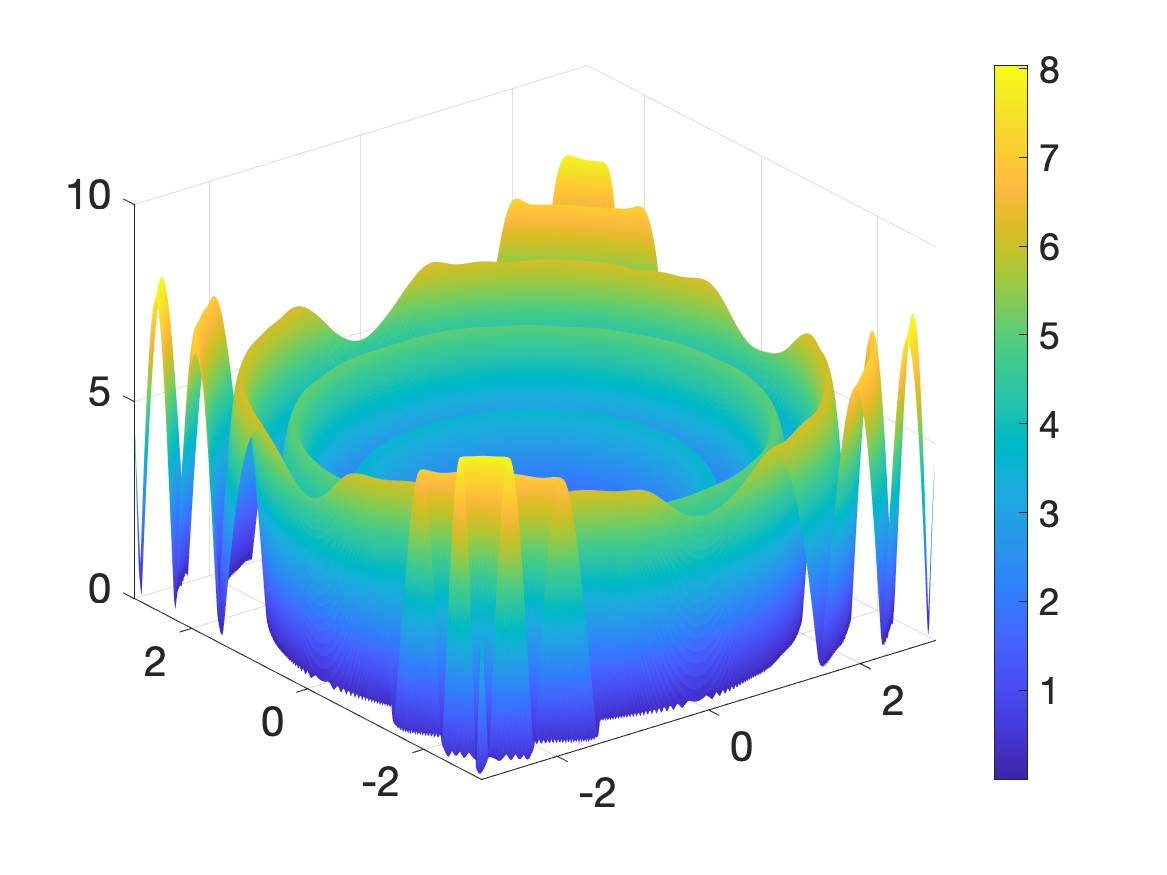}}
	\quad
	\subfloat[\label{4c}$|\nabla f^*- \nabla f_{\rm comp}|/\|\nabla f^*\|_{L^{\infty}}$]{\includegraphics[width = .3\textwidth]{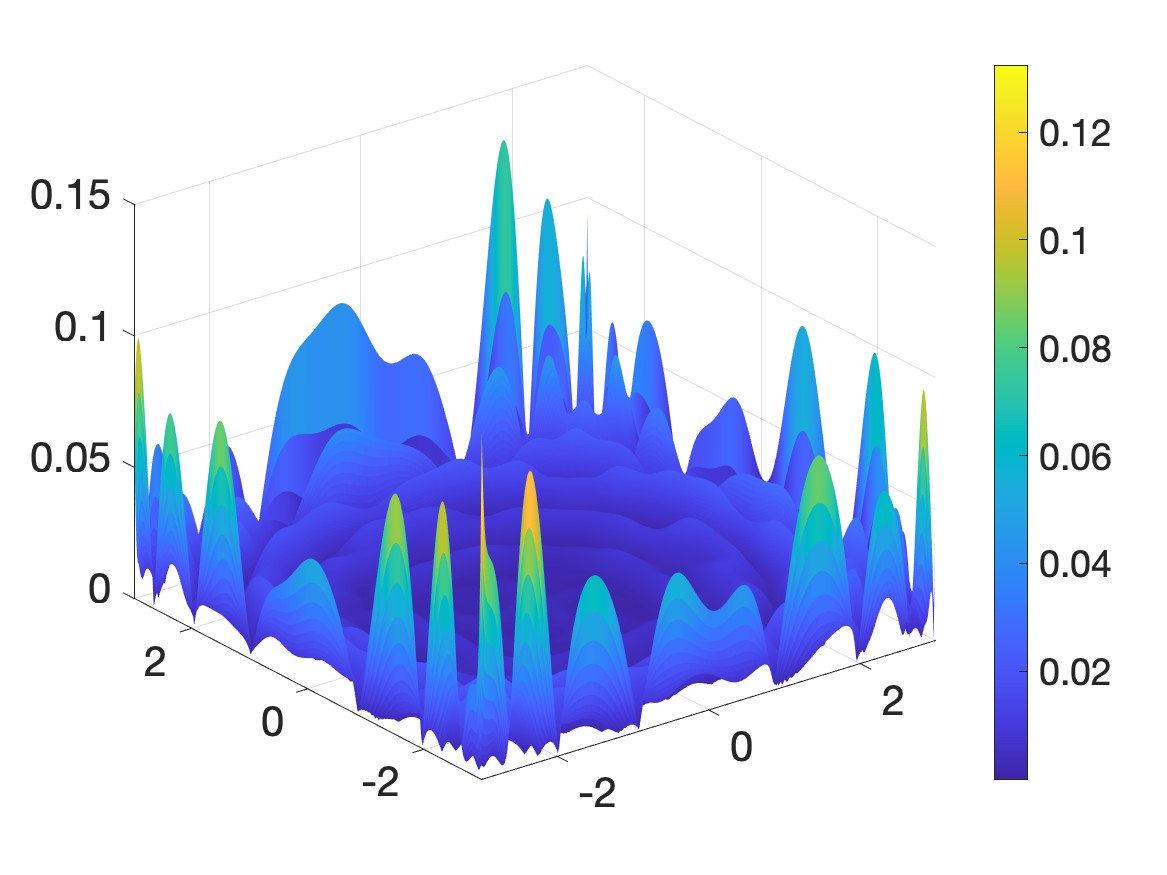}}
	
	\subfloat[The true function $\Delta f^*$]{\includegraphics[width = .3\textwidth]{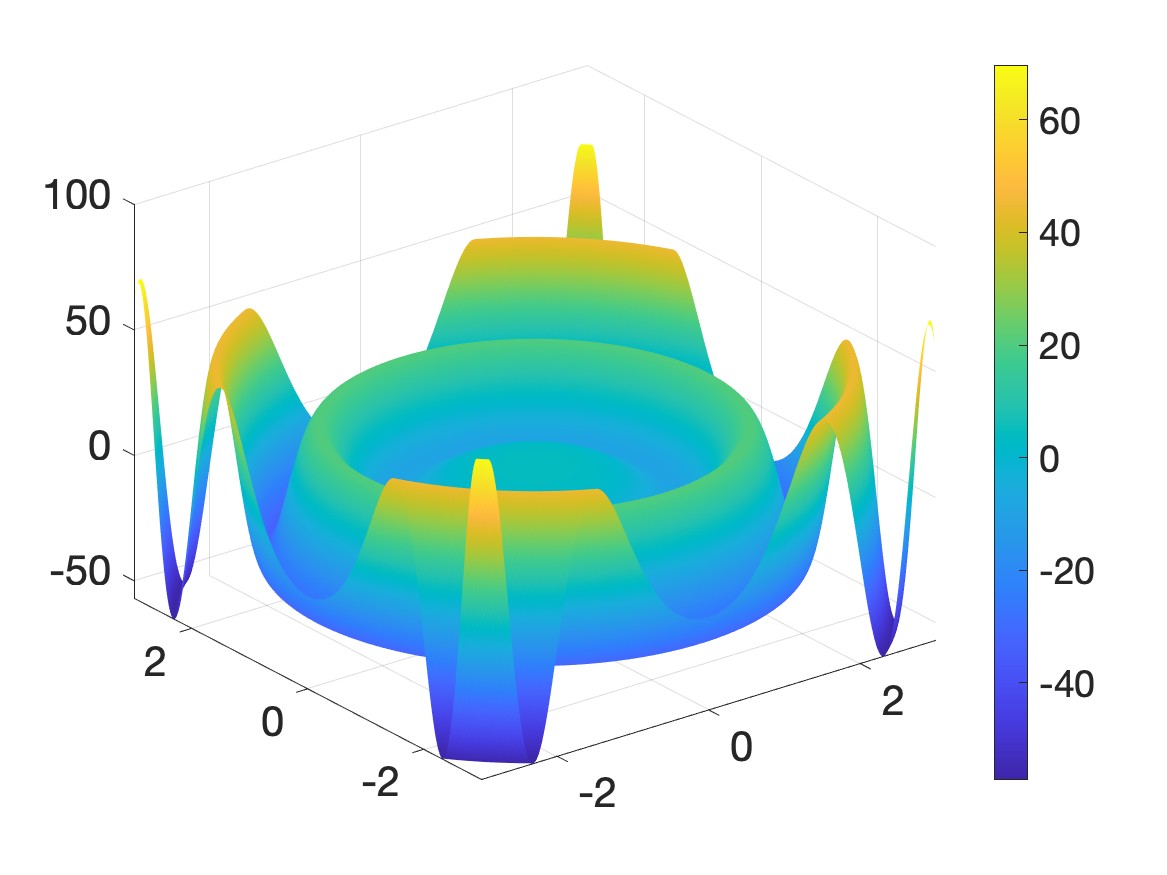}}
	\quad
	\subfloat[The computed function $\Delta f_{\rm comp}$]{\includegraphics[width = .3\textwidth]{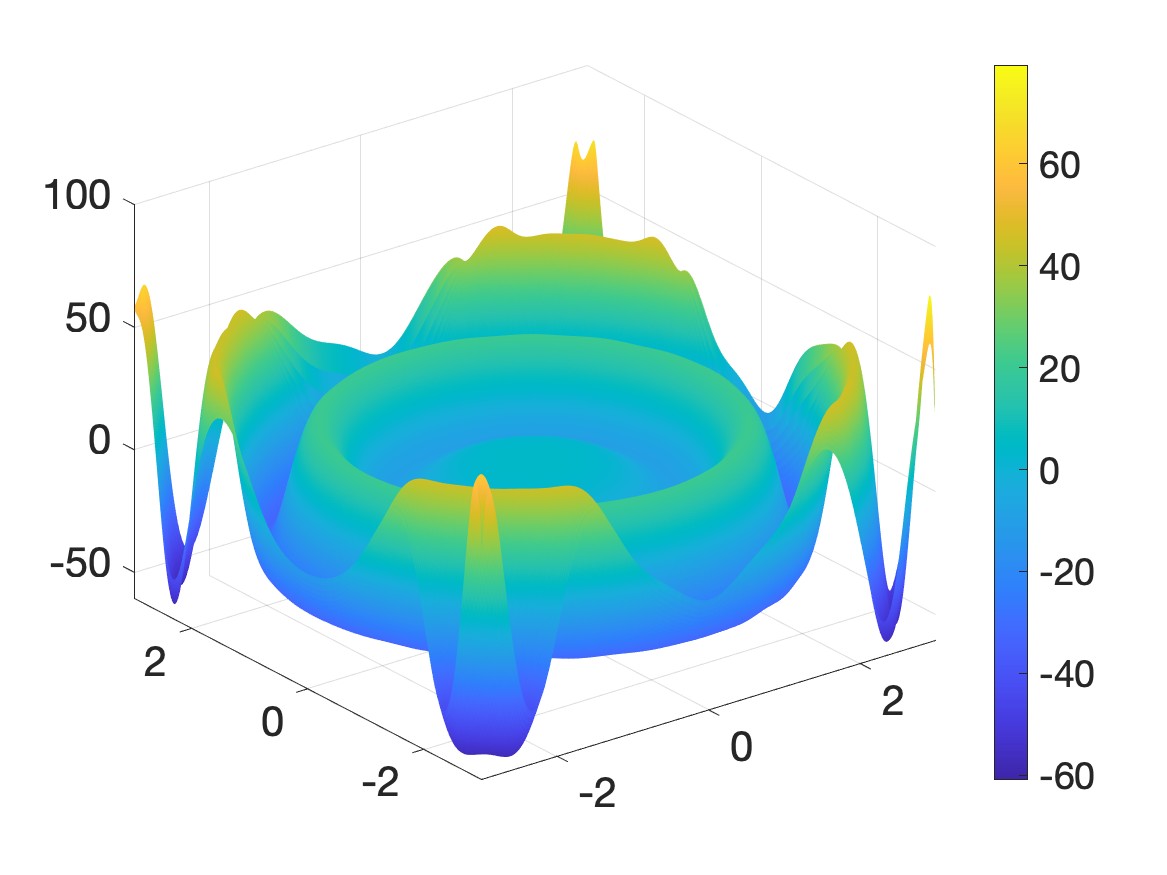}}
	\quad
	\subfloat[\label{4f}$|\Delta f^*- \Delta f_{\rm comp}|/\|\Delta f^*\|_{L^{\infty}}$]{\includegraphics[width = .3\textwidth]{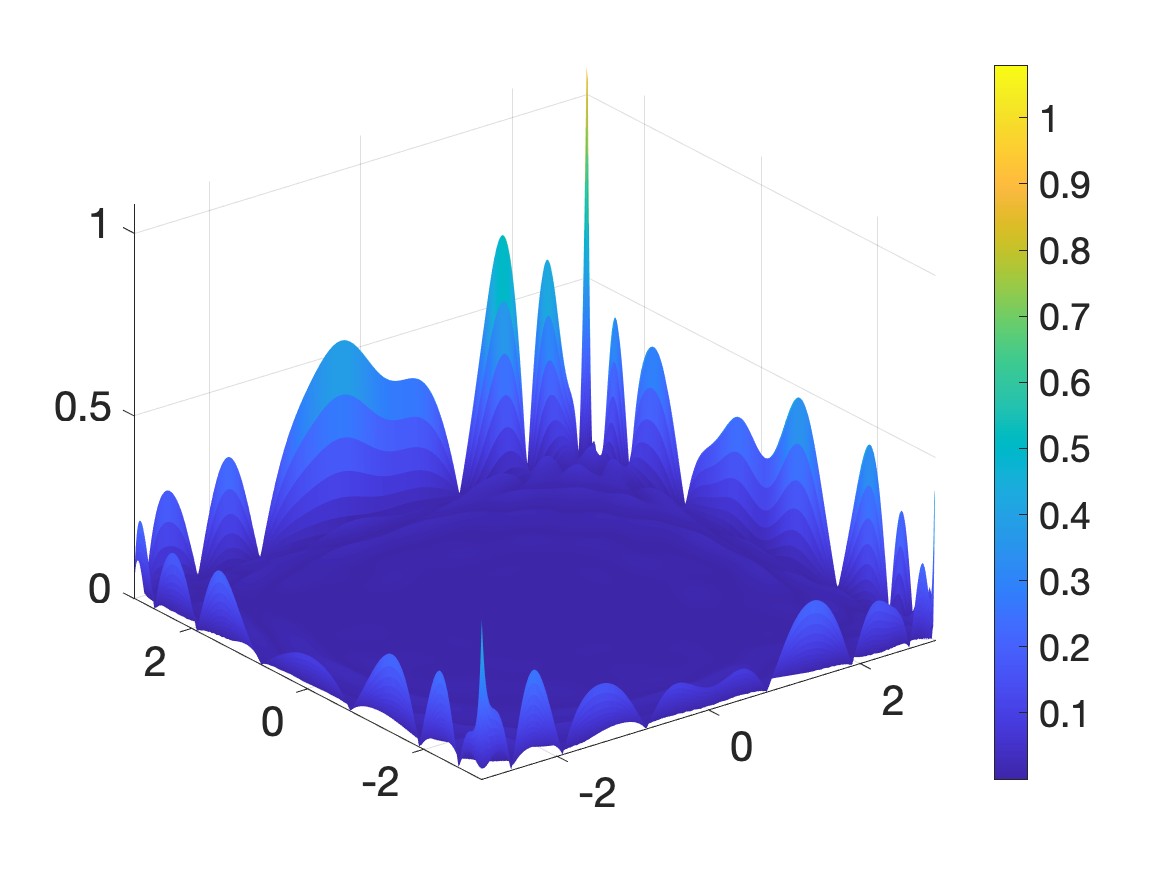}}
	
	\caption{\label{fig2d1} Test 3. The true and computed $|\nabla f|$ and $\Delta f$,  using data that has been corrupted with a 10\% noise level.}
\end{figure}

In this test, we take $f^*(x, y) = \sin(x^2 + y^2)$. Then
\begin{align*}
	|\nabla f^*(x, y)| &= 2\sqrt{(x^2 + y^2) \cos^2(x^2 + y^2)},\\
	\Delta f^*(x,y) &= 4\cos(x^2 + y^2) -4x^2\sin(x^2 + y^2) - 4y^2\sin(x^2 + y^2).  
\end{align*}
In Figure \ref{fig2d1}, we show the true and computed magnitudes of the gradient, $|\nabla f^*|$, and the Laplacian, $\Delta f^*$, are depicted using data that has been affected by a 10\% noise level.

It follows from Figure \ref{fig2d1} that we successfully compute the first and second derivatives of the function $f^*(x, y) = \sin(x^2 + y^2).$ It is interesting to mention that the errors in the computation are less than the noise level.
As seen in the last column of Figures \ref{4c} and \ref{4f}, the computed errors occur mainly on the boundary of the computed domain.
This feature is true when we test our algorithm with up to the noise level 20\%.

\begin{table}[ht]
    \centering
\begin{tabular}{|c|c|c|c|c|}
    \hline 
    		&$\frac{\|\nabla f_{\rm true} - \nabla f_{\rm comp}\|_{L^2(\Omega)}}{\|\nabla f_{\rm true}\|_{L^2(\Omega)}}$
    		& $\frac{\|\Delta f_{\rm true} - \Delta f_{\rm comp}\|_{L^2(\Omega)}}{\|\Delta f_{\rm true}\|_{L^2(\Omega)}}$
		&$\frac{\|\nabla f_{\rm true} - \nabla f_{\rm comp}\|_{L^2(\Omega')}}{\|\nabla f_{\rm true}\|_{L^2(\Omega')}}$
		&  $\frac{\|\Delta f_{\rm true} - \Delta f_{\rm comp}\|_{L^2(\Omega')}}{\|\Delta f_{\rm true}\|_{L^2(\Omega')}}$
  \\  \hline
    		$\delta = 0.05$
		&	 0.0303	
		& 	0.1282
		& 0.0163
		&0.0320
  \\  \hline
    		$\delta = 0.10$
		&0.0313
		&0.1306
		&0.0165
		&0.0324
	  \\  \hline
    		$\delta = 0.20$
		& 0.0338
		&0.1610
		&0.0173
		& 0.0337
	\\\hline
\end{tabular}
\caption{\label{tab3} Test 3. The relative computed errors with respect to the $L^2$ norms on $\Omega = (-3, 3)^2$ and $ \Omega' = (-2, 2)^2$. As can be observed, the relative errors are lower than the noise level.}
\end{table}

The performance of our method with different noise levels $\delta \in \{5\%, 10\%, 20\%\}$ is out of expectation. Table \ref{tab3} shows that the relative computed errors are consistent with the corresponding noise levels. The relative error with respect to the $L^2(\Omega')$ norm is significantly low.

\subsubsection{Test 2}

\begin{figure}[ht]
	\subfloat[The true function $|\nabla f^*|$]{\includegraphics[width = .3\textwidth]{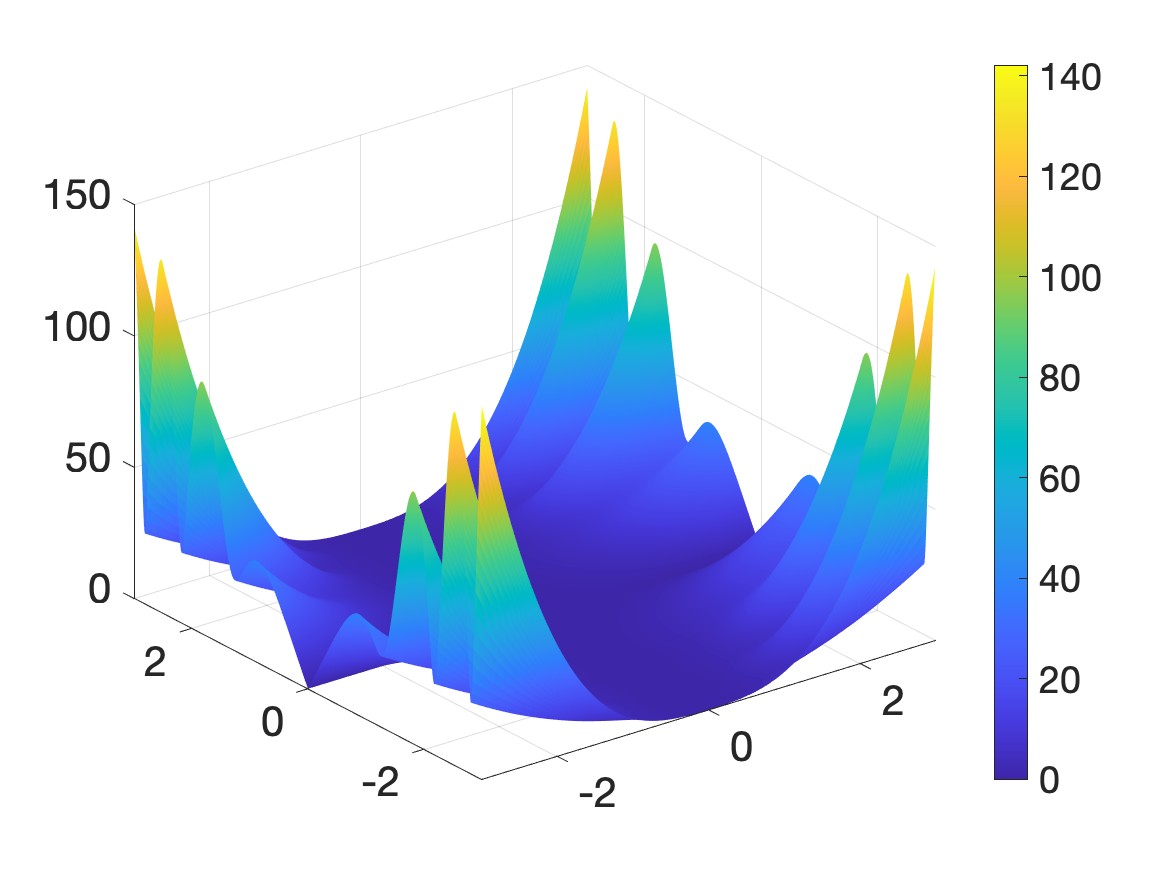}}
	\quad
	\subfloat[The computed function $|\nabla f_{\rm comp}|$]{\includegraphics[width = .3\textwidth]{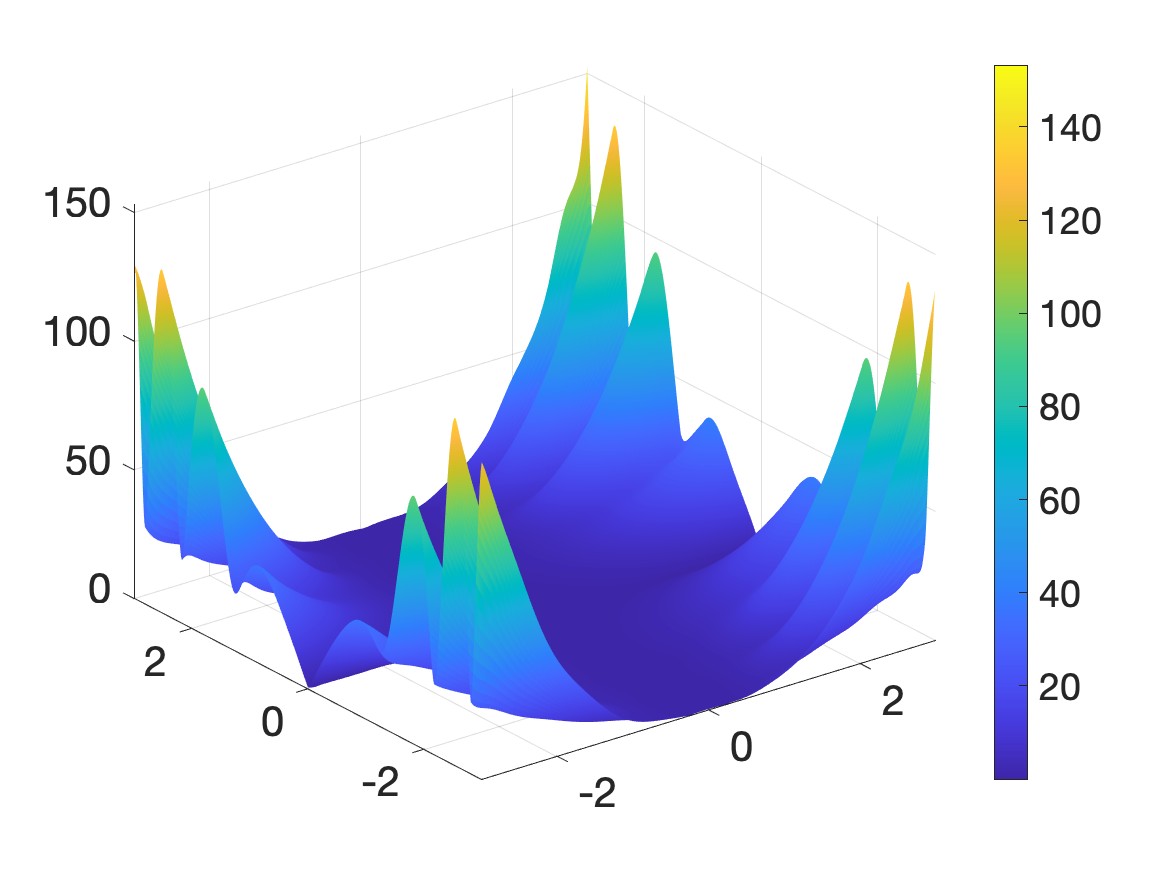}}
	\quad
	\subfloat[$|\nabla f^*- \nabla f_{\rm comp}|/\|\nabla f^*\|_{L^{\infty}}$]{\includegraphics[width = .3\textwidth]{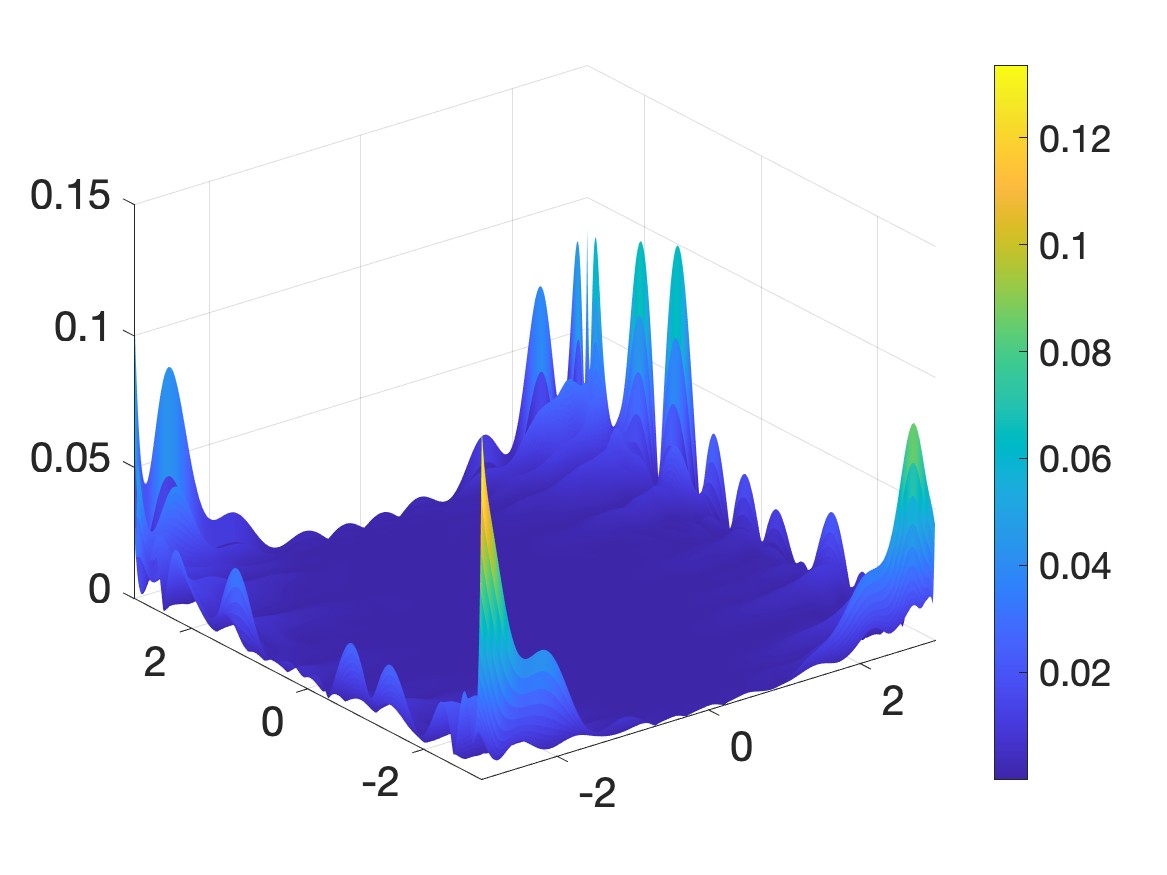}}
	
	\subfloat[The true function $\Delta f^*$]{\includegraphics[width = .3\textwidth]{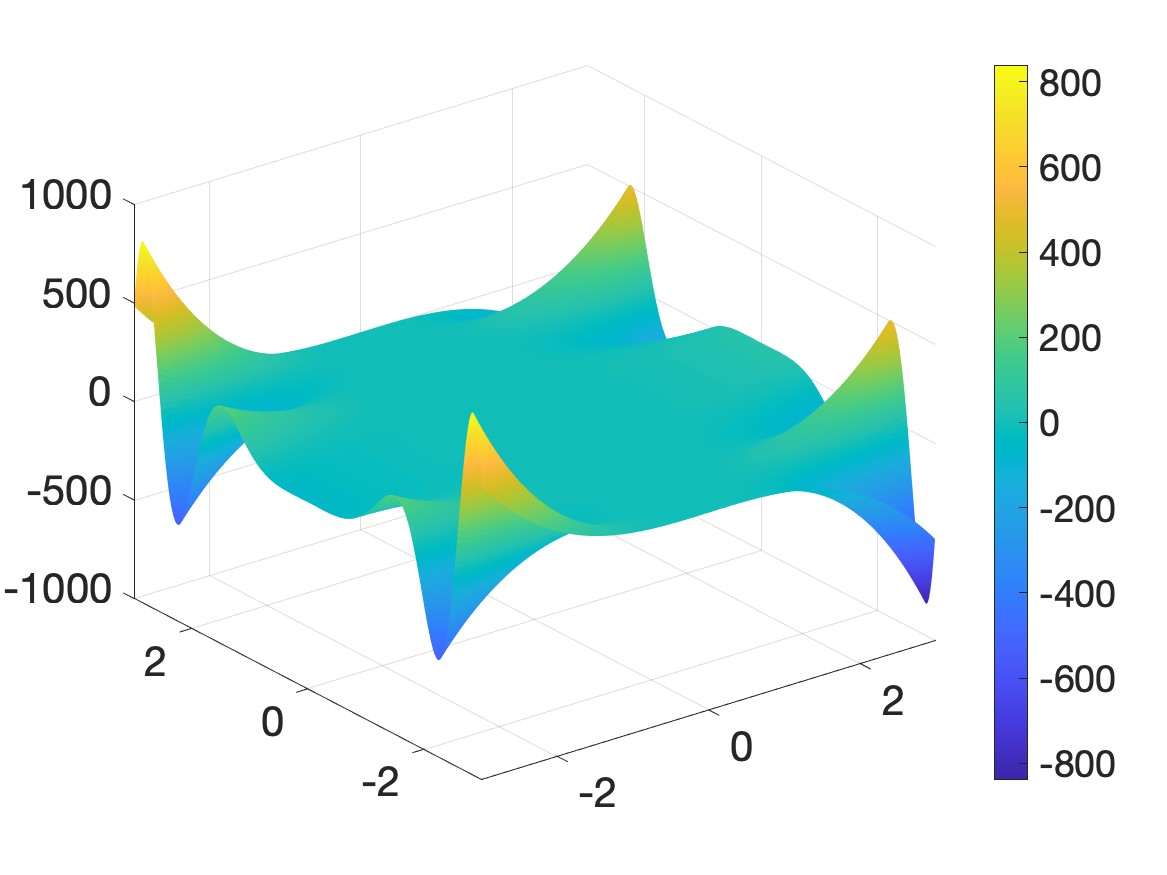}}
	\quad
	\subfloat[The computed function $\Delta f_{\rm comp}$]{\includegraphics[width = .3\textwidth]{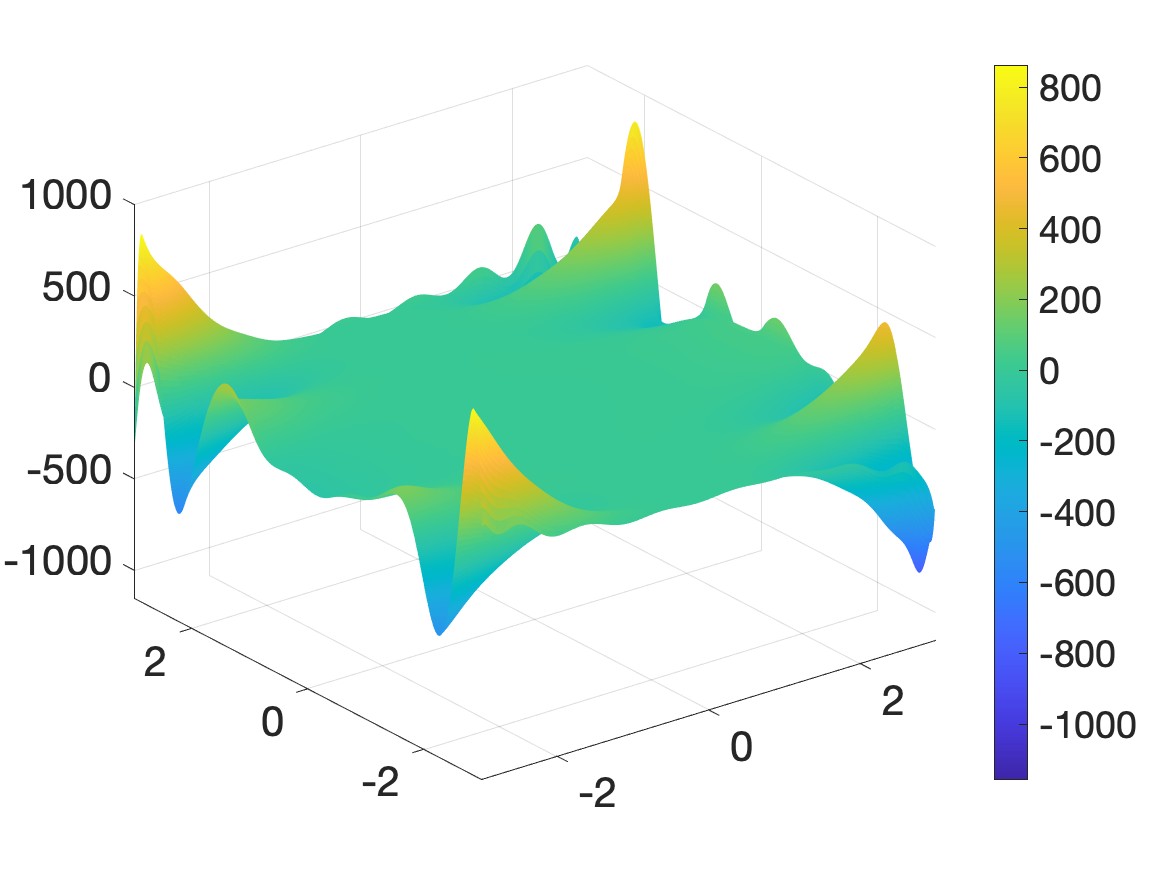}}
	\quad
	\subfloat[$|\Delta f^*- \Delta f_{\rm comp}|/\|\Delta f^*\|_{L^{\infty}}$]{\includegraphics[width = .3\textwidth]{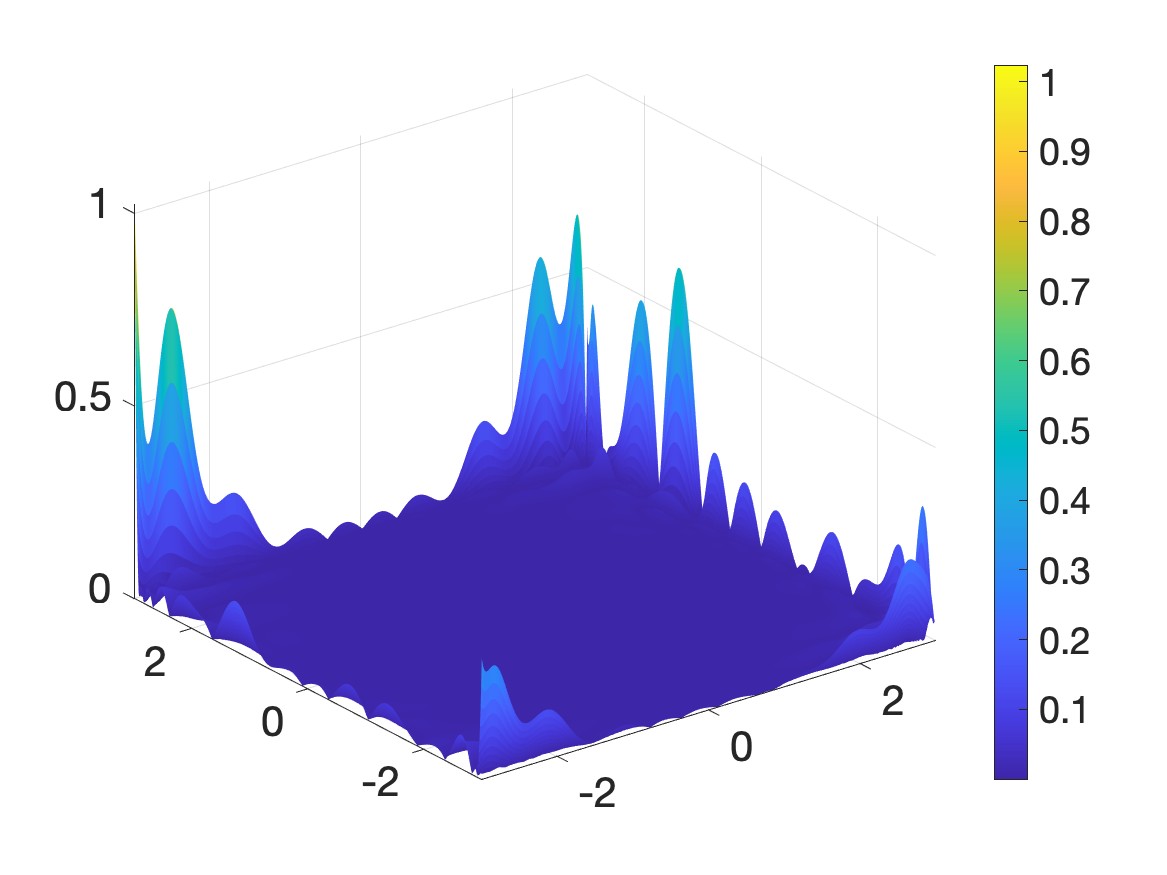}}
	
	\caption{\label{fig2d2} Test 4. The true and computed $|\nabla f|$ and $\Delta f$,  using data that has been corrupted with a 10\% noise level.}
\end{figure}

We test our algorithms by computing derivatives of $f^*(x, y) = x^3 \sin(y^2).$
As in Test 2.1, we will compute numerical versions of the functions
\begin{align*}
	|\nabla f^*(x, y)| &= \sqrt{9x^4 \sin^2(y^2) + 4x^6y^2\cos^2(y^2)}
	\\
	\Delta f^*(x, y) &= (-4x^3 y^2 + 6x) \sin(y^2) + 2x^3 \cos(y^2).
\end{align*}
The numerical results are displayed in Figure \ref{fig2d2}.

Like Test 3, it is evident from Figure \ref{fig2d2} and Table \ref{tab4} that Algorithm \ref{alg} provides out-of-expectation numerical results for the first and second derivatives of the function $f^*$ with acceptable relative errors.

\begin{table}[ht]
    \centering
\begin{tabular}{|c|c|c|c|c|}
    \hline 
    		&$\frac{\|\nabla f_{\rm true} - \nabla f_{\rm comp}\|_{L^2(\Omega)}}{\|\nabla f_{\rm true}\|_{L^2(\Omega)}}$
    		& $\frac{\|\Delta f_{\rm true} - \Delta f_{\rm comp}\|_{L^2(\Omega)}}{\|\Delta f_{\rm true}\|_{L^2(\Omega)}}$
		&$\frac{\|\nabla f_{\rm true} - \nabla f_{\rm comp}\|_{L^2(\Omega')}}{\|\nabla f_{\rm true}\|_{L^2(\Omega')}}$
		&  $\frac{\|\Delta f_{\rm true} - \Delta f_{\rm comp}\|_{L^2(\Omega')}}{\|\Delta f_{\rm true}\|_{L^2(\Omega')}}$
  \\  \hline
    		$\delta = 0.05$
		&0.0336	 	
		& 0.1981	
		& 0.0128
		&0.0587
  \\  \hline
    		$\delta = 0.10$
		&0.0386
		&0.2301
		&0.0142
		&0.0645
	  \\  \hline
    		$\delta = 0.20$
		&0.0571
		&0.3666
		&0.0154
		&0.0734
	\\\hline
\end{tabular}
\caption{\label{tab4} Test 4. The relative computed errors with respect to the $L^2$ norms on $\Omega = (-3, 3)^2$ and $ \Omega' = (-2, 2)^2$. As can be observed, the relative errors are lower than the noise level.}
\end{table}

\section{Comparison with some existing differentiating methods}\label{sec5}

In this section, we compare the efficiency of Algorithm \ref{alg} with some other widely-used methods for differentiation.
We brevity, we only implement and present the computation of these methods in one dimension.

\subsection{Fourier expansion with the trigonometric basis}\label{sec4.1}

{\it If we substitute the polynomial-exponential basis with the trigonometric basis in Algorithm \ref{alg}, what would be the outcome?}
 Let $N_{\rm trig} > 0$ be a cut-off number. 
Write
\begin{equation}
	f^{\delta}(x) \simeq f_{\rm trig}^{\delta}(x) =  a_0 + \sum_{n = 1}^{N_{\rm trig}} a_n \cos\Big(\frac{\pi n x}{R}\Big) + \sum_{n = 1}^{N_{\rm trig}} b_n \sin \Big(\frac{n \pi x}{R}\Big)
	\label{Fourier}
\end{equation}
where 
\begin{align*}
	a_0 &= \frac{1}{2R}  \int_{-R}^R f^{\delta}(s)ds,\\
	a_n &= \frac{1}{R} \int_{-R}^R f^{\delta}(s) \cos\Big(\frac{\pi n s}{R}\Big)ds, \quad n \geq 1,\\
	b_n &= \frac{1}{R} \int_{-R}^R f^{\delta}(s) \sin\Big(\frac{\pi n s}{R}\Big)ds, \quad n \geq 1.
\end{align*}
Like in Section \ref{sec3.1}, we choose the optimal cut-off number $N_{\rm trig}$ using the following formula
\begin{equation}
	N_{\rm trig} = \underset{3 \leq N \leq 25}{\rm argmin}
	 \Big\{
		\|f^{\delta} - \Big[a_0 + \sum_{n = 1}^{N} a_n \cos\Big(\frac{\pi n x}{R}\Big) + \sum_{n = 1}^{N} b_n \sin \Big(\frac{n \pi x}{R}\Big)\Big]\|_{L^{\infty}(-R, R)}
	\Big\}.
	\label{4.2}
\end{equation}
The range $3 \leq N \leq 25$ in \eqref{4.2} can be changed.
However, we observe that when $N$ exceeds 25, the numerical results get worse.
Using \eqref{Fourier}, we approximate ${f^*}'$ via the formula
\begin{equation}
	f_{\rm trig}'(x) \simeq - \sum_{n = 1}^{N_{\rm trig}}  \frac{\pi n a_n}{R} \sin\Big(\frac{\pi n x}{R}\Big) + \sum_{n = 1}^{N_{\rm trig}} \frac{n \pi b_n}{R} \cos \Big(\frac{n \pi x}{R}\Big).
	\label{5.2}
\end{equation}
We also can approximate ${f^*}''$ via
\begin{equation}
	f_{\rm trig}''(x)  =  -\sum_{n = 1}^{N_{\rm trig}} \Big(\frac{\pi n}{R}\Big)^2 a_n \cos\Big(\frac{\pi n x}{R}\Big) - \sum_{n = 1}^{N_{\rm trig}} \Big(\frac{\pi n}{R}\Big)^2 b_n \sin \Big(\frac{n \pi x}{R}\Big).
	\label{5.3}
\end{equation}

%

In rows 3 of Table \ref{tab6} and Table \ref{tab7}, we present the relative errors of the method utilizing \eqref{5.2} and \eqref{5.3} with respect to the $L^2(-3, 3)$ norm. The experiment employs data with 5\%, 10\%, and 20\% noise components. The comparison of these results and those via Algorithm \ref{alg} in rows 2 and 5 reveals that, in terms of accuracy, the polynomial-exponential basis outperforms the trigonometric basis. We conclude that the use of the polynomial-exponential basis is crucial to the success of our approach.

\subsection{Comparison with the Tikhonov regularization method}\label{sec4.2}

\begin{figure}[ht]
\begin{center}
	\subfloat[The graph of the function $l(\gamma)$, $\gamma \in (10^{-5}, 10^{-2})$, computed by 5\% noisy data taken in Test 1]{\includegraphics[width = .3\textwidth]{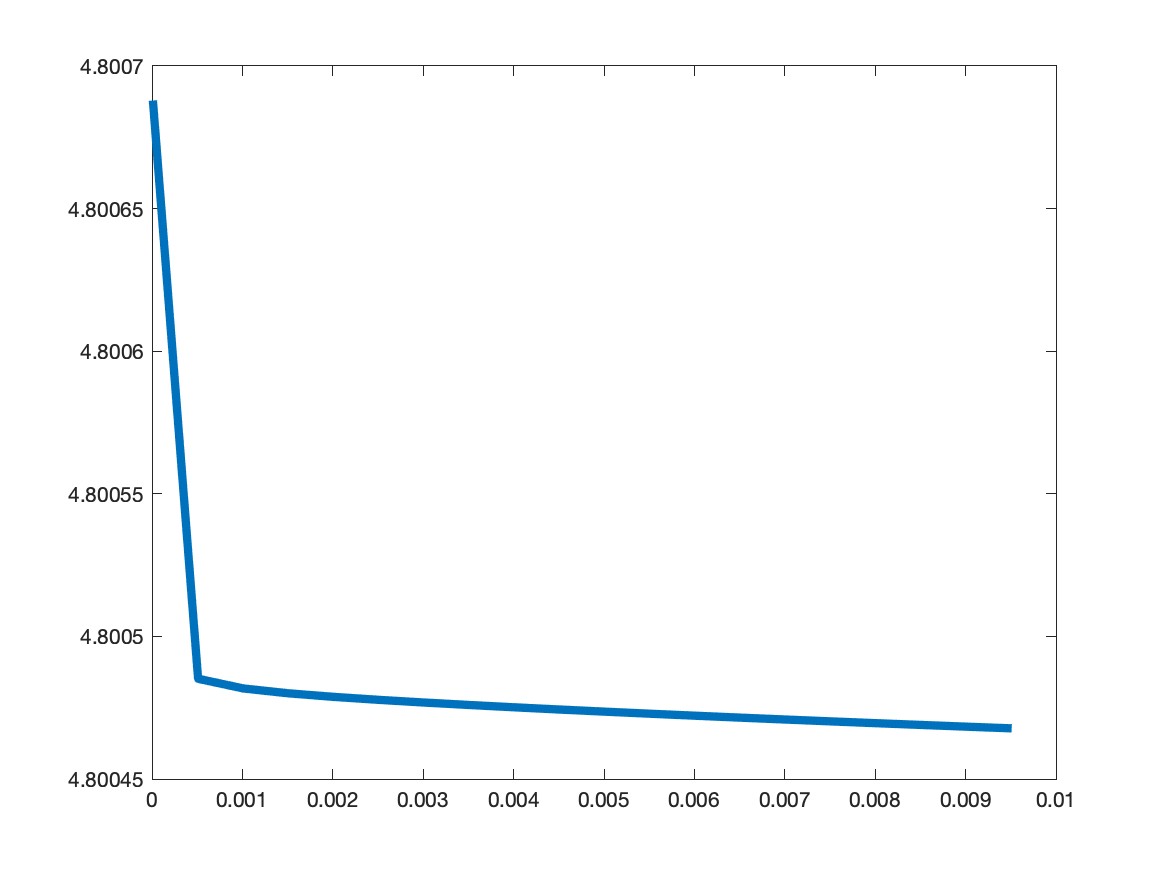}}
	\quad
	\subfloat[The graph of the function $l(\gamma)$,  $\gamma \in (10^{-5}, 10^{-2})$, computed by 5\% noisy data taken in Test 2]{\includegraphics[width = .3\textwidth]{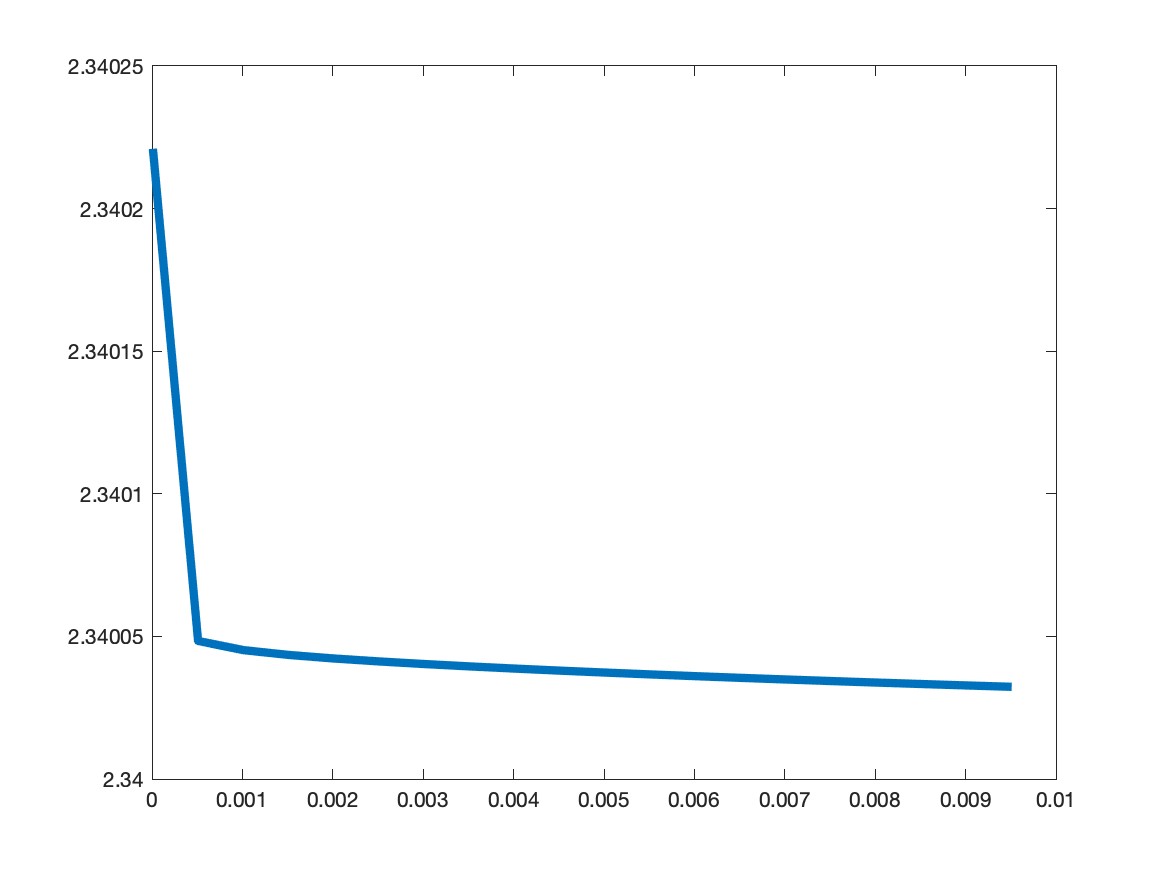}}
\end{center}
\caption{\label{Lshape}The $L$-shape behavior of the function $l$. The optimal regularization $\gamma$ occurs at the corner of the graph. }
\end{figure}

A widely used approach to differentiate the data is to employ the Tikhonov optimization. Define $u = f'$.
By the fundamental theorem of calculus, we write 
\[
	f^{\delta}(x) - f^{\delta}(-R) = \int_{-R}^x u(s)ds.
\]
This suggests us to find $u$ by minimizing the following Tikhonov cost functional
\[
	J_\gamma(u) = \int_{-R}^R \Big|\int_{-R}^x u(s)ds - [f^{\delta}(x) - f^{\delta}(-R)]\Big|^2dx
	+ \gamma \|u\|_{L^2(-R, R)}^2
\]
where $\gamma > 0$ is a regularization parameter. The term $\gamma \|u\|_{L^2(-R, R)}^2$ is called the regularization term. 
It is possible to employ higher-order norms, such as $H^1$ and $H^2$, instead of the $L^2$ norm for the regularization term. The derivative ${f^{\delta}_{\rm Tik}}'$ is determined to be $u^{\gamma}_{\rm min}$, which is the minimizer of $J_\gamma$. 
We compute the second derivative of $f$, denoted by ${f^{\delta}_{\rm Tik}}''$, by differentiating ${f^{\delta}_{\rm Tik}}'$ using the Tikhonov optimization approach again.
While the Tikhonov regularization approach is known for its robustness, selecting the appropriate norm and regularization parameter is not trivial.
 A popular approach to selecting the regularization parameter is Morosov's discrepancy principle \cite{Morozov:springer1984}, in which
the regularization parameter $\gamma$ is determined such that
\begin{equation}
	 \Big\|\int_{-R}^x u^{\gamma}_{\rm min}(s) ds - [f(x) - f(-R)]\Big\|_{L^2(-R, R)} = O(\delta)
\end{equation}
where $\delta$ is the noise level.
Morosov's discrepancy method requires knowledge of the noise level, which is not always known in practical applications. 
The $L-$shape method is a good option to find $\gamma$. One sketches the graph of the function $l: (0, 1) \to \R$ defined as 
\[
	l(\gamma) =  \Big\|\int_{-R}^x u^{\gamma}_{\rm min}(s) ds - [f(x) - f(-R)]\Big\|_{L^2(-R, R)}^2.
\]
The graph of $l$ is an $L$-shape curve, see Figure \ref{Lshape}.
The optimal regularization parameter is determined by the horizontal coordinate of the corner of the $L$ curve.
We show in row 4 and row 8 of Table \ref{tab6} the relative computing error using the Tikhonov regularization method, in which we use the $L$-shape method to choose the parameter $\gamma.$
Comparing the relative errors in row 2 and row 6 of this Table, we conclude that our method performs better than the Tikhonov optimization method.

\begin{remark}
	Since we need to compute $l(\gamma)$, we have to repeatedly attempt to minimize $J_{\gamma}$ for each value of $\gamma$, which makes the optimization-based method time-consuming.
In the conducted experiments, we utilized an iMac 3.2 GHz Quad-Core Intel Core I5 to evaluate $l(\gamma)$ on uniformly partitioning the interval $(10^{-5}, 10^{-2})$ with a step size of $5\times10^{-4}$. The computation took 93.85 seconds.
In contrast, the same computer required just 0.01167 seconds to compute both the first and second derivatives of the function in each of Test 1 and Test 2.
Thus, our method is considerably faster.
\end{remark}

\subsection{Comparison with the method based on cubic spline}\label{spline}

One can apply cubic spline curves, introduced in \cite{Reinsch:nm2967, Reinsch:nm1971, Schoenberg:PAMS1964}, to smooth data before differentiating it. 
Define a partition $\{x_i = -R + (i - 1)h, i = 1, 2, \dots, N\} $ of the computational interval $(-R, R)$ similarly to the partition in \eqref{3.2222}.
The difference between this partition and the one in \eqref{3.2222} is that we change the step size to $h = 0.1$ because the cubic spline method works better with sparse data.
Then, we use the "spline" command in Matlab to construct a set of spline curves that fit the data $f^{\delta}$.
The command spline provides the coefficients of each cubic function in each subdomain $[x_i, x_{i + 1}]$, $i \in \{1, \dots, N_x-1\}$, that allows us to compute the first and second derivatives of the function $f^{\delta}.$

To evaluate the performance of the cubic spline method, we show the relative error in computation in rows 5  of Table \ref{tab6} and Table \ref{tab7}. Our experiment reveals that this method is well-suited for computing the first derivative; however, it may lead to a large computation error when computing the second derivative. A comparison of the first and final columns of the aforementioned tables leads us to the conclusion that our method yields markedly more accurate results.

\subsection{Numerical results}

\begin{figure}[ht]
    \centering
    \subfloat[Test 1; noise level $\delta = 5\%$]{\includegraphics[width = .3\textwidth]{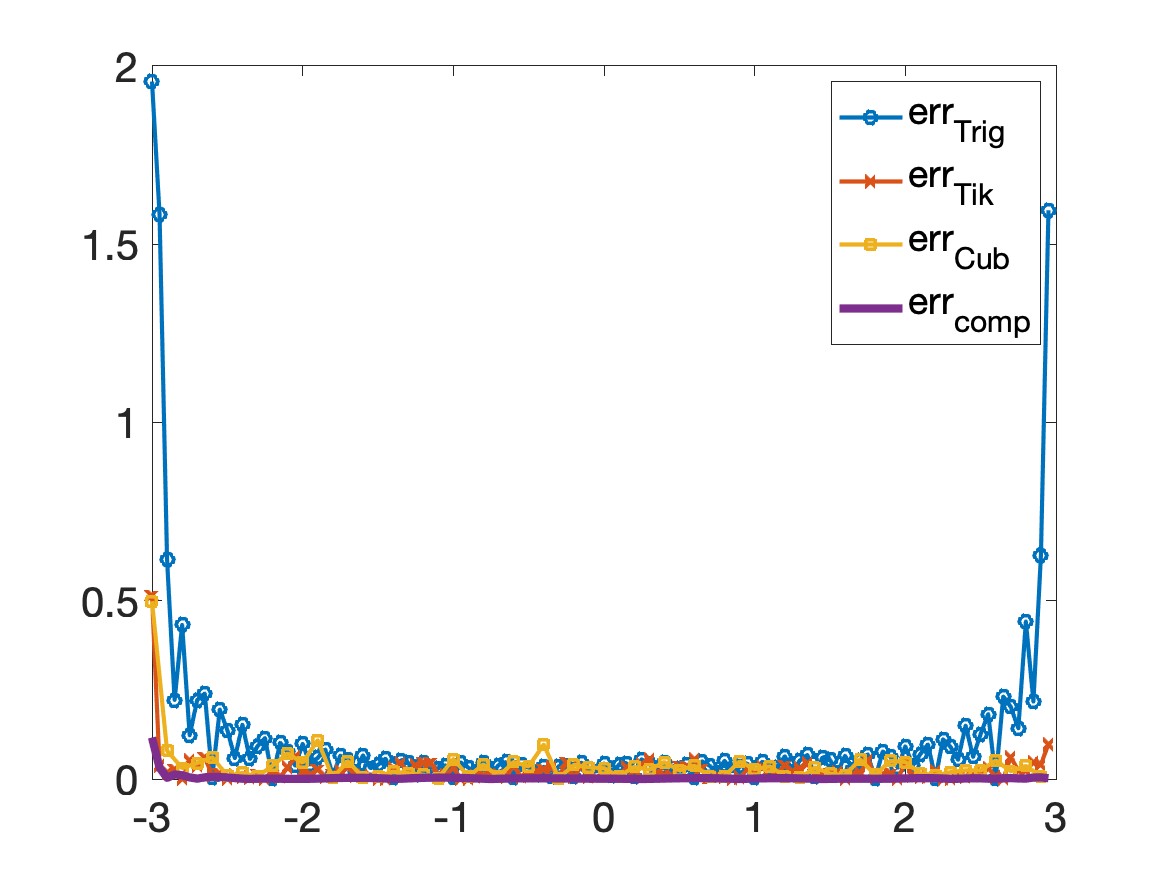}}
    \quad
    \subfloat[Test 1; noise level $\delta = 10\%$]{\includegraphics[width = .3\textwidth]{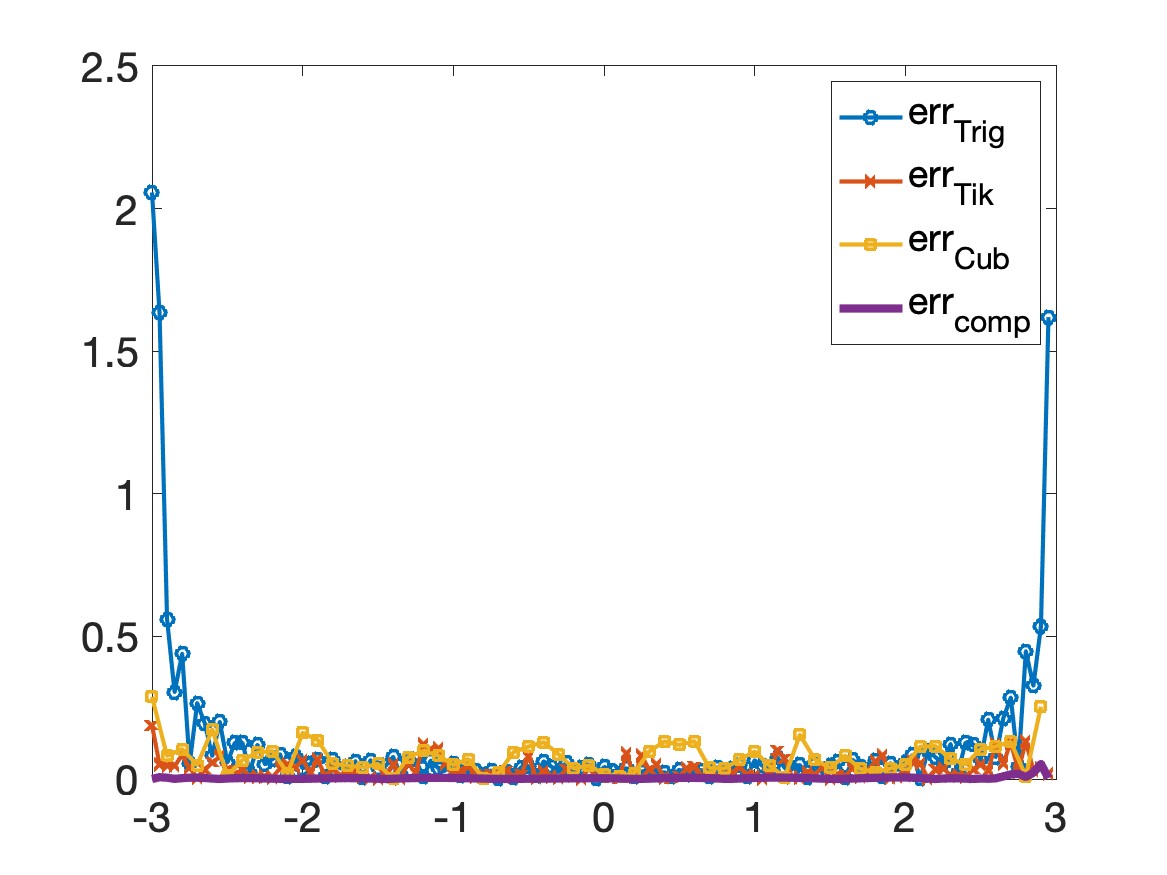}}
    \quad
    \subfloat[Test 1; noise level $\delta = 20\%$]{\includegraphics[width = .3\textwidth]{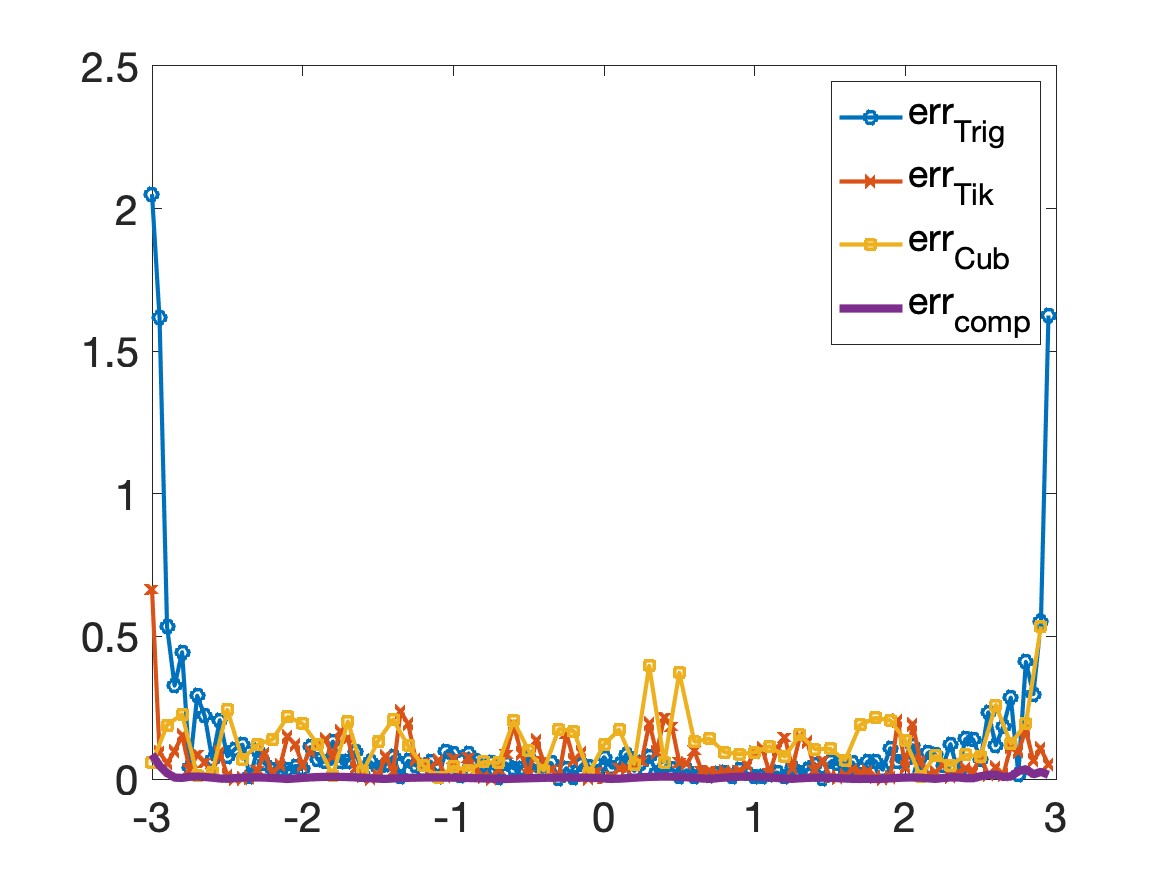}}

\subfloat[Test 2; noise level $\delta = 5\%$]{\includegraphics[width = .3\textwidth]{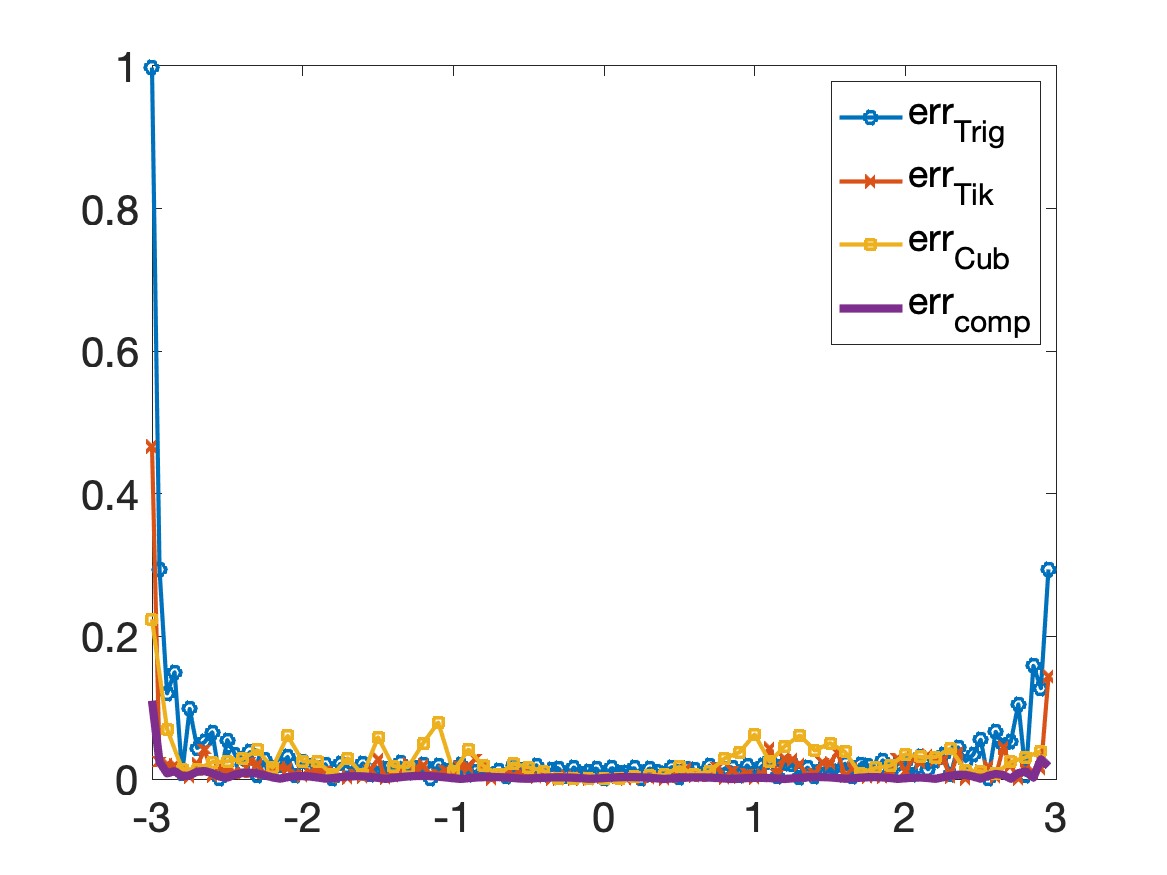}}
    \quad
    \subfloat[Test 2; noise level $\delta = 10\%$]{\includegraphics[width = .3\textwidth]{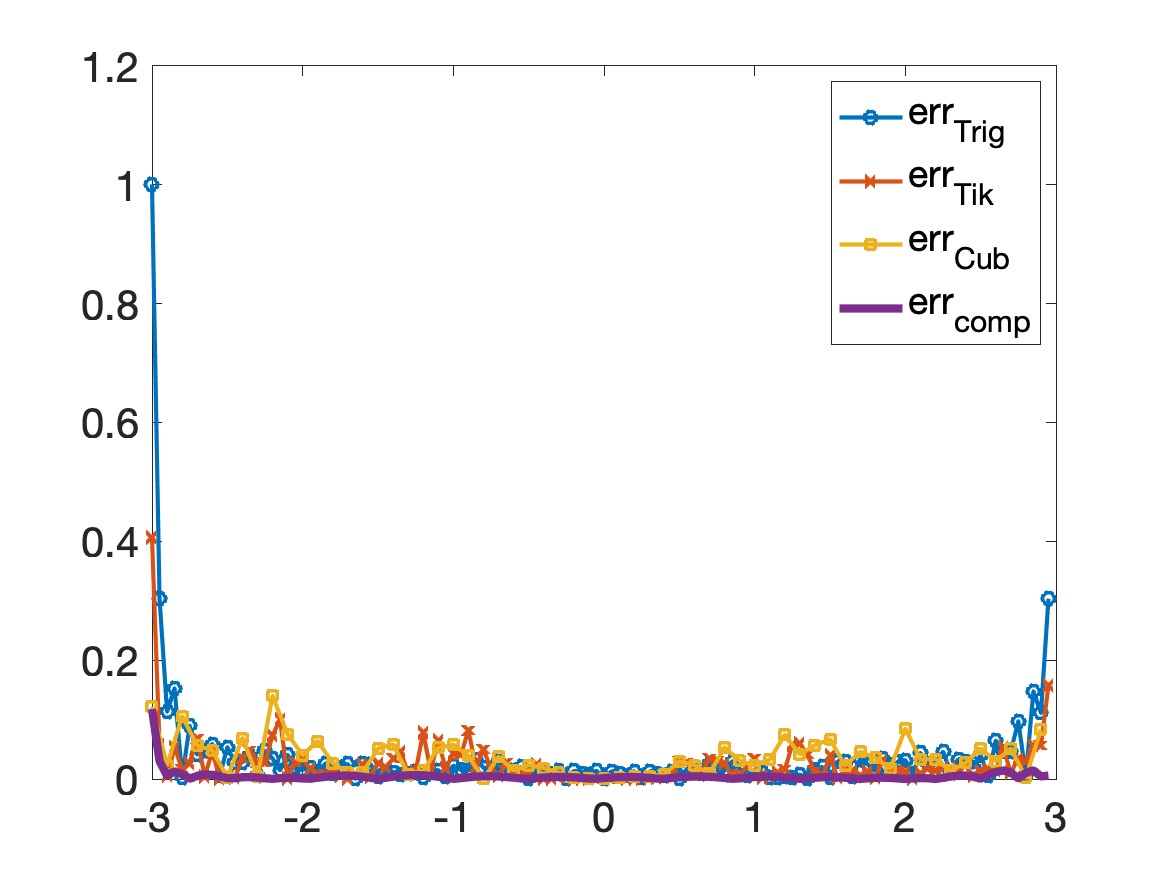}}
    \quad
    \subfloat[Test 2; noise level $\delta = 20\%$]{\includegraphics[width = .3\textwidth]{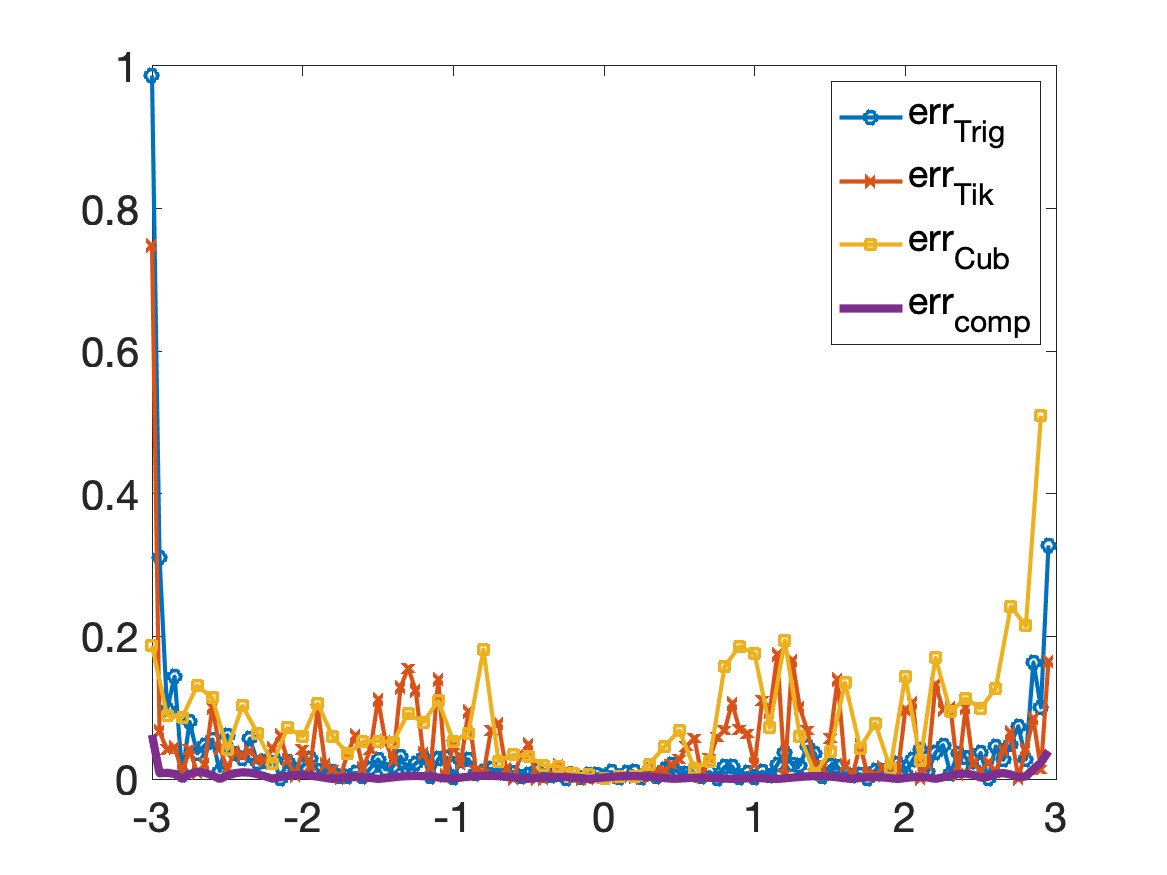}}
    
    \caption{\label{fig7} The relative difference of the true and computed first derivatives of the functions in Test 1 and Test 2. These relative difference functions are defined in \eqref{4.6} for $x \in (-3, 3)$.}   
\end{figure}

In Figure \ref{fig7}, we show the graphs of the following functions 
\begin{equation}
	\begin{array}{cc}
		\ds{\rm err}_{\rm Trig}(x) = \frac{|f_{\rm Trig}'(x) - f_{\rm true}'(x)|}{\|f_{\rm true}'\|_{L^{\infty}(-3, 3)}},  
	&\ds {\rm err}_{\rm Tik}(x) = \frac{|f_{\rm Tik}'(x) - f_{\rm true}'(x)|}{\|f_{\rm true}'\|_{L^{\infty}(-3, 3)}},\\
	\ds{\rm err}_{\rm Cub}(x) = \frac{|f_{\rm Cub}'(x) - f_{\rm true}'(x)|}{\|f_{\rm true}'\|_{L^{\infty}(-3, 3)}},  
	& \ds{\rm err}_{\rm comp}(x) = \frac{|f_{\rm comp}'(x) - f_{\rm true}'(x)|}{\|f_{\rm true}'\|_{L^{\infty}(-3, 3)}}
	\end{array}
	\label{4.6}
\end{equation}
where $f_{\rm Trig}'$, $f_{\rm Tik}'$, $f_{\rm Cub}'$, and $f_{\rm comp}'$ are the computed first derivatives via the methods in Sections \ref{sec4.1}, \ref{sec4.2}, \ref{spline}, and \ref{section_comp}, respectively.
The data is taken from Test 1 and Test 2 in Section \ref{sec4.1} when the noise level are 5\%, 10\%, and 20\%.
 We do not show similar graphs for the computed derivatives since the methods mentions in Sections \ref{sec4.1}, \ref{sec4.2}, and \ref{spline} are not accurate, see Table \ref{tab6}. 
It is evident from Figure \ref{fig7} that our method performs better than the methods mentioned.

\begin{table}[ht]
    \centering
\begin{tabular}{|c|c|c|c|c|c|c|}
\hline
&Test 1
&Test 1
&Test 1 
&Test 2
&Test 2
&Test 2
\\
&
$\delta = 0.05$
&$\delta = 0.10$
&$\delta = 0.20$
&$\delta = 0.05$
&$\delta = 0.10$
&$\delta = 0.20$
\\
\hline
 $\frac{\|f'_{\rm true} - f'_{\rm comp}\|_{L^2(-3, 3)}}{\|f'_{\rm true}\|_{L^2(-3, 3)}}$
&{\bf 0.0060}
&{\bf 0.0110}
&{\bf 0.0260}
&{\bf 0.0017}
&{\bf 0.0047}
&{\bf 0.0117}
\\
\hline
$\frac{\|f'_{\rm true} - f'_{\rm trig}\|_{L^2(-3, 3)}}{\|f'_{\rm true}\|_{L^2(-3, 3)}}$
&0.4157
&0.4397
&0.4485
&0.2216
&0.2219
&0.2260
\\
\hline
$\frac{\|f'_{\rm true} - f'_{\rm Tik}\|_{L^2(-3, 3)}}{\|f'_{\rm true}\|_{L^2(-3, 3)}}$
&0.0652
&0.0879
&0.1484
&0.1115
&0.1248
&0.1780
\\
\hline
$\frac{\|f'_{\rm true} - f'_{\rm cub}\|_{L^2(-3, 3)}}{\|f'_{\rm true}\|_{L^2(-3, 3)}}$
&0.0775
&0.1185
&0.2450
&0.0551
&0.1904
&0.2355
\\
\hline
   \end{tabular}
\caption{\label{tab6} The relative errors in computing the first derivatives of our method and several other approaches are evaluated with respect to the $L^2(-3, 3)$ norm. In this table, $f'_{\rm true}$ is the true first derivative that we aim to approximate. The function $f'_{\rm comp}$ is the first derivative computed by our approach. The functions $f'_{\rm trig}$, $f'_{\rm Tik}$ and $f'_{\rm cub}$ are the derivative computed by the existing methods reviewed in Sections \ref{sec4.1}, \ref{sec4.2}, and \ref{spline} respectively.}
\end{table}

\begin{table}[ht]
    \centering
\begin{tabular}{|c|c|c|c|c|c|c|}
\hline
&Test 1
&Test 1
&Test 1 
&Test 2
&Test 2
&Test 2
\\
&
$\delta = 0.05$
&$\delta = 0.10$
&$\delta = 0.20$
&$\delta = 0.05$
&$\delta = 0.10$
&$\delta = 0.20$
\\
\hline
$\frac{\|f''_{\rm true} - f''_{\rm comp}\|_{L^2(-3, 3)}}{\|f''_{\rm true}\|_{L^2(-3, 3)}}$
&{\bf 0.0268}
&{\bf 0.0996}
&{\bf 0.1123}
&{\bf 0.0380}
&{\bf 0.0955}
&{\bf 0.1734}
\\
\hline
$\frac{\|f''_{\rm true} - f''_{\rm trig}\|_{L^2(-3, 3)}}{\|f''_{\rm true}\|_{L^2(-3, 3)}}$
&1.2419
&1.4411
&1.5455
&1.0224
&1.0314
&1.0674
\\
\hline
$\frac{\|f''_{\rm true} - f''_{\rm Tik}\|_{L^2(-3, 3)}}{\|f''_{\rm true}\|_{L^2(-3, 3)}}$
&0.3410
&0.4504
&0.6334
&0.5912
&0.6044
&0.6700
\\
\hline
$\frac{\|f'_{\rm true} - f''_{\rm cub}\|_{L^2(-3, 3)}}{\|f''_{\rm true}\|_{L^2(-3, 3)}}$
&0.5406
&0.6970
&1.2402
&0.5118
&0.7455
&1.0360
\\
\hline
   \end{tabular}
\caption{\label{tab7} The relative errors in computing the second derivatives of our method and several other approaches are evaluated with respect to the $L^2(-3, 3)$ norm. In this table, $f''_{\rm true}$ is the true second derivative that we aim to approximate. The function $f''_{\rm comp}$ is the second derivative computed by our approach. The functions $f''_{\rm trig}$, $f''_{\rm Tik}$ and $f''_{\rm cub}$ are the derivative computed by the existing methods reviewed in Sections \ref{sec4.1}, \ref{sec4.2}, and \ref{spline} respectively.}
\end{table}

We also show in Table \ref{tab6} and Table \ref{tab7} the relative computational errors with respect to the $L^2(-3,3)$ norm. 
	Based on the information in Table \ref{tab6} and Table \ref{tab7}, we can infer that our approach is the most precise one. The relative errors in computing the first derivative are roughly ten times lower than the noise levels, while the other methods have relative errors that are consistent with the noise levels. Additionally, our method's relative errors in computing second derivatives are in line with the noise levels, whereas other methods' relative errors in computing second derivatives are not satisfactory.


\section{Concluding remarks}\label{sec6}

We introduce a technique to compute derivatives of data affected by noise, which can be difficult due to the potential for significant computation errors even from a small amount of noise. Our method is straightforward. Firstly, we represent the data by its Fourier series with respect to the polynomial-exponential basis. Then, we remove all high-frequency components by truncating this series. The desired derivatives are then obtained directly. To demonstrate the efficacy of our algorithm, we provide numerical examples in both one and two dimensions. Furthermore, we compare our approach with some commonly used methods and conclude that our method is more precise.


\begin{thebibliography}{10}

\bibitem{Ahnert:cpc2007}
K.~Ahnert and L.~A. Segel.
\newblock Numerical differentiation of experimental data: local versus global
  methods.
\newblock {\em Computer Physics Communications}, 177:764--774, 2007.

\bibitem{Breugel:IEEEAccess}
F.~V. Breugel, J.~N. Kutz, and B.~W. Brunton.
\newblock Numerical differentiation of noisy data: {A} unifying multi-objective
  optimization framework.
\newblock {\em IEEE Access}, 8:196865--196877, 2020.

\bibitem{Engl:Kluwer1996}
H.~W. Engl, M.~Hanke, and A.~Neubauer.
\newblock {\em Regularization of Inverse Problems, Mathematics and its
  Applications}.
\newblock Kluwer Academic Publishers Group, Dordrecht, 1996.

\bibitem{Friedrichs:tams1944}
K.~O. Friedrichs.
\newblock The identity of weak and strong extensions of differential operators.
\newblock {\em Trans. Amer. Math. Soc.}, 55:132--151, 1944.

\bibitem{Groetsch:amm1991}
C.~W. Groetsch.
\newblock Differentiation of approximately specified functions.
\newblock {\em Am. Math. Mon.}, 98:847--850, 1991.

\bibitem{Hanke2001}
M.~Hanke and O.~Sherzer.
\newblock Inverse problems light numerical differentiation.
\newblock {\em Am. Math. Mon.}, 108:512--521, 2001.

\bibitem{VoKlibanovNguyen:IP2020}
V.~A. Khoa, G.~W. Bidney, M.~V. Klibanov, L.~H. Nguyen, L.~Nguyen, A.~Sullivan,
  and V.~N. Astratov.
\newblock Convexification and experimental data for a {3D} inverse scattering
  problem with the moving point source.
\newblock {\em Inverse Problems}, 36:085007, 2020.

\bibitem{Khoaelal:IPSE2021}
V.~A. Khoa, G.~W. Bidney, M.~V. Klibanov, L.~H. Nguyen, L.~Nguyen, A.~Sullivan,
  and V.~N. Astratov.
\newblock An inverse problem of a simultaneous reconstruction of the dielectric
  constant and conductivity from experimental backscattering data.
\newblock {\em Inverse Problems in Science and Engineering}, 29(5):712--735,
  2021.

\bibitem{KhoaKlibanovLoc:SIAMImaging2020}
V.~A. Khoa, M.~V. Klibanov, and L.~H. Nguyen.
\newblock Convexification for a 3{D} inverse scattering problem with the moving
  point source.
\newblock {\em SIAM J. Imaging Sci.}, 13(2):871--904, 2020.

\bibitem{Klibanov:jiip2017}
M.~V. Klibanov.
\newblock Convexification of restricted {D}irichlet to {N}eumann map.
\newblock {\em J. Inverse and Ill-Posed Problems}, 25(5):669--685, 2017.

\bibitem{KlibanovLeNguyen:SIAM2020}
M.~V. Klibanov, T.~T. Le, and L.~H. Nguyen.
\newblock Convergent numerical method for a linearized travel time tomography
  problem with incomplete data.
\newblock {\em SIAM Journal on Scientific Computing}, 42:B1173--B1192, 2020.

\bibitem{KlibanovLiBook}
M.~V. Klibanov and J.~Li.
\newblock {\em Inverse Problems and Carleman Estimates: Global Uniqueness,
  Global Convergence and Experimental Data}.
\newblock De Gruyter, 2021.

\bibitem{KlibanovNguyen:ip2019}
M.~V. Klibanov and L.~H. Nguyen.
\newblock {PDE}-based numerical method for a limited angle {X}-ray tomography.
\newblock {\em Inverse Problems}, 35:045009, 2019.

\bibitem{KlibanovTimonov:cac2023}
M.~V. Klibanov and A.~Timonov.
\newblock A comparative study of two globally convergent numerical methods for
  acoustic tomography.
\newblock {\em Communications in Analysis and Computation}, 1:12--31, 2023.

\bibitem{Knowles:ejde2014}
I.~Knowles and R.~J. Renka.
\newblock Methods for numerical differentiation of noisy data.
\newblock {\em Elect. J. Diff. Equations}, Conference 21:235--246, 2014.

\bibitem{Knowles:nm1995}
I.~Knowles and R.~Wallace.
\newblock A variational method for numerical differentiation.
\newblock {\em Numer. Math.}, 70:91--110, 1995.

\bibitem{LeNguyen:JSC2022}
T.~T. Le and L.~H. Nguyen.
\newblock The gradient descent method for the convexification to solve boundary
  value problems of quasi-linear {PDEs} and a coefficient inverse problem.
\newblock {\em Journal of Scientific Computing}, 91(3):74, 2022.

\bibitem{Morozov:springer1984}
V.~A. Mozorov.
\newblock {\em Methods for solving incorrectly posed problems}.
\newblock Springer Verlag, New York, 1984.

\bibitem{Ramlau:jnfao2002}
R.~Ramlau.
\newblock Morozov's discrepancy principle for {T}ikhonov regularization of
  nonlinear operators.
\newblock {\em Journal for Numer. Funct. Anal. and Opt.}, 23:147--172, 2002.

\bibitem{Ramm:mc2001}
A.~G. Ramm and A.~B. Smirnova.
\newblock On stable numerical differentiation.
\newblock {\em Math. Comp.}, 70:1131--1153, 2001.

\bibitem{Reinsch:nm2967}
C.~H. Reinsch.
\newblock Smoothing by spline functions.
\newblock {\em Numer. Math.}, 10:177--183, 1967.

\bibitem{Reinsch:nm1971}
C.~H. Reinsch.
\newblock Smoothing by spline functions. {II}.
\newblock {\em Numer. Math.}, 16:451--454, 1971.

\bibitem{Scherzer:SIAMJna1993}
O.~Scherzer.
\newblock The use of {M}orozov's discrepancy principle for {T}ikhonov
  regularization for solving non-linear ill-posed problems.
\newblock {\em SIAM J. Numer. Anal.}, 30:1796--1838, 1993.

\bibitem{Schoenberg:PAMS1964}
I.~J. Schoenberg.
\newblock Spline functions and the problem of graduation.
\newblock {\em Proc. Amer. Math. Soc.}, 52:497--950, 1964.

\end{thebibliography}
\end{document}